\documentclass[12pt]{article}
\usepackage{amsmath,amssymb,amsfonts,amscd,amsxtra}
\usepackage{latexsym}

\input{epsf.sty}
\usepackage{epsfig}
\usepackage{graphicx,color}
\usepackage{graphics}
\usepackage{subfigure}
\usepackage{appendix}
\usepackage[top=1in, bottom=1in, left=1in, right=1in]{geometry}

\renewcommand{\theequation}{\thesection.\arabic{equation}}

\newtheorem{thm}{Theorem}[section]

\newtheorem{lem}[thm]{Lemma}
\newtheorem{prop}[thm]{Proposition}

\newtheorem{rmk}[thm]{Remark}

\renewcommand{\Im}{{\mbox{Im}}}
\renewcommand{\Re}{{\mbox{Re}}}

\newcommand{\abs}[1]{\left\vert#1\right\vert}

\newcommand{\qed}{\hfill \ensuremath{\square}}

\renewcommand\appendix{\par
  \setcounter{section}{0}
  \setcounter{subsection}{0}
  \setcounter{figure}{0}
  \setcounter{table}{0}
  \renewcommand\thesection{Appendix \Alph{section}}
  \renewcommand\theequation{\Alph{section}.\arabic{equation}}
  \renewcommand\thefigure{\Alph{section}.\arabic{figure}}
  \renewcommand\thetable{\Alph{section}.\arabic{table}}
  \renewcommand\thethm{\Alph{section}.\arabic{thm}}
}

\numberwithin{equation}{section}

\date{}

\title{Scattering by a periodic array of subwavelength slits I: field enhancement in the diffraction regime}
\author{
Junshan Lin \thanks{\footnotesize Department of Mathematics and Statistics, Auburn University, Auburn, AL 36849 (jzl0097@
auburn.edu). Junshan Lin was partially supported by the NSF grant DMS-1417676.}
 \; and Hai Zhang\thanks{\footnotesize 
Department of Mathematics, 
 HKUST,  Clear Water Bay, Kowloon, Hong Kong (haizhang@ust.hk). Hai Zhang was supported by HK RGC grant ECS 26301016 and the UGC grant SBI17SC12 from HKUST.}}

\begin{document}

\maketitle

\begin{abstract}
This is the first part in a series of two papers that concern with the quantitative analysis of the electromagnetic field enhancement and anomalous diffraction by a periodic array of subwavelength slits. The scattering problem in the diffraction regime is investigated in this part, for which the size of the period is comparable to the incident wavelength. We distinguish scattering resonances and real eigenvalues, and derive their asymptotic expansions when they are away from the Rayleigh cut-off frequencies. Furthermore, we present quantitative analysis of the field enhancement at resonant frequencies, by quantifying both the enhancement order and the associated resonant modes.
 The field enhancement near the Rayleigh cut-off frequencies is also investigated. It is demonstrated that the field enhancement
 becomes weaker at the resonant frequency if it is close to the Rayleigh cut-off frequencies. Finally, we also characterize the embedded eigenvalues 
 for the underlying periodic structure, and point out that transmission anomaly such as Fano resonant phenomenon does not occur for the narrow slit array.
 \end{abstract}


\setcounter{equation}{0}
\setlength{\arraycolsep}{0.25em}
\section{Introduction}\label{sec-introduction}
There has been increasing interest in the electromagnetic scattering by subwavelength apertures or holes in recent years,
due to their significant applications in biological and chemical sensing, near-field spectroscopy, etc 
\cite{abajo, astilean00, chen14, ebbesen98, garcia10, krie04, lin14, lin15,seo09, sturman10, taka01, yang02}. 
Such subwavelenth structures can generate extraordinary optical transmission and strongly enhanced local electromagnetic fields,
which leads extremely high sensitivity in the sensing or imaging of biological or chemical samples.
However, as of today there are still controversies over the mechanisms contributing to the anomalous field enhancement  \cite{garcia10}.
The complication has to do with the multiscale nature of the underlying metallic structures as well as various enhancement behaviors that it induces.
For instance, the enhancement can be attributed to surface plasmonic resonance \cite{ebbesen98, garcia10}, non-plasmonic resonances \cite{taka01, yang02}, or even without the resonant effect (cf. \cite{lin14, lin15,lin_zhang16})

Very recently, for a single narrow slit perforated in a perfect conducting slab, we have presented quantitative analysis of the field enhancement,
which provides a complete picture for its enhancement mechanisms \cite{lin_zhang16}. In this paper and its sequel \cite{lin_zhang17}, 
we investigate the scattering and
field enhancement when the slab is patterned with a periodic array of narrow slits. The physics becomes 
richer compared to the single slit case. In addition, the quantitative studies of field enhancement also present new mathematical challenges.
(i) While the scattering problem for the single slit attains a unique solution, there exists a set of singular frequencies (real eigenvalues)
for the periodic strucuture;  (ii) In contrast to the single slit case, the field enhancement factor for the periodic case would depend on the size of period and the incident angle for a given frequency; (iii) The scattering by the periodic structure will exhibit so-called Rayleigh anomaly when certain propagating spatial harmonics turns into an evanescent mode, which gives rise to diffractive anomaly near such cut-off frequencies; (iv) 
The field enhancement mechanisms in the diffraction regime and the homogenization regime are different and need separate investigations.
The incident wavelength is comparable to the size of  the period in the former case, but much larger than the size of the period in the latter.
As to be demonstrated in this paper and  \cite{lin_zhang17}, while the field enhancement in the diffraction regime is mostly attributed to the resonant phenomenon, the enhancement and diffraction anomaly in the homogenization regimes can be induced by non-resonant phenomenon.

\begin{figure}[!htbp]
\begin{center}
\includegraphics[height=5.5cm,width=15cm]{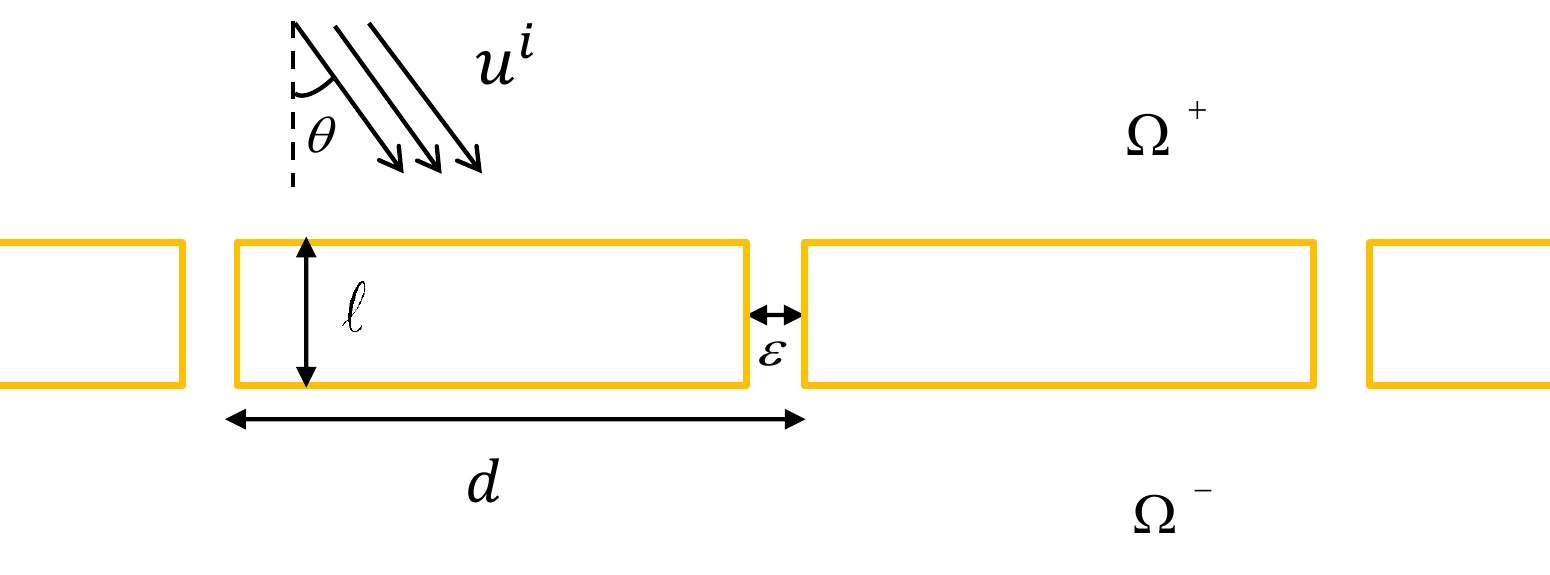}
\caption{Geometry of the scattering problem. 
The slits  $S_\varepsilon$ are arranged periodically with the size of the period $d$, and each slit has a rectangular shape of length $\ell$ and width $\varepsilon$ respectively. 
The domains above and below the perfect conductor slab are denoted as $\Omega^{+}$ and $\Omega^{-}$ respectively, and the domain exterior to the perfect conductor is denoted as $\Omega_{\varepsilon}$, which consists of $S_\varepsilon$, $\Omega^{+}$, and $\Omega^{-}$. }\label{fig-prob_geo}
\end{center}
\end{figure}

We now present the setup of our scattering problem. We consider a perfect conducting slab that is perforated with a periodic array of slits
and the geometry of its cross section is depicted in Figure \ref{fig-prob_geo}. The slab occupies the domain 
$\{(x_1,x_2)\;|\; 0<x_2<\ell\}$ on the $x_1x_2$ plane, where $\ell$ is the thickness of the metal. The slits, which are invariant along the $x_3$ direction, 
occupy the region $\displaystyle{S_\varepsilon=\bigcup_{n=0}^{\infty} S_\varepsilon^{(0)} + nd}$, where $d$ is the size of the period, and
$S_\varepsilon^{(0)}:=\{(x_1,x_2)\;|\; 0<x_1<\varepsilon, 0<x_2<\ell \}$ has a rectangular shape.
Let us denote the semi-infinite domain above and below the slab by $\Omega^{+}$ and $\Omega^{-}$.
We also denote by $\Omega_{\varepsilon}$ be the domain exterior to  the perfect conductor, i.e.,
$\Omega_{\varepsilon}=\Omega^{+}\cup\Omega^{-}\cup S_\varepsilon$, and $\nu$ denotes the unit outward normal pointing to the exterior domain $\Omega^{+}$ or $\Omega^{-}$.

We are interested in the case when the width of the slit $\varepsilon$ is much smaller than the thickness of the slab $\ell$.
In this paper, we shall assume that the size of period $d$ has the same order as the thickness of the slab $\ell$.
Furthermore, the incident wavelength $\lambda$ is comparable to $d$ such that  the problem under consideration is in the diffraction regime. 
Namely, $\varepsilon \ll \ell$ and $d  \sim \ell \sim \lambda$ hold for the scaling of parameters for the problem under consideration.
For clarity, we shall assume $\ell=1$ in all technical derivations throughout the paper. 
The enhancement theory for the case of $\ell\neq1$ then follows by a normalization process and a scaling argument.

Assume that a polarized time-harmonic electromagnetic wave impinges upon the perfect conductor from the above.
We consider the transverse magnetic (TM) case where the incident magnetic field is perpendicular to the  $x_1x_2$ plane,
and its $x_3$ component is given by the scalar function $u^{i}= e^{i k(\sin \theta x_1 - \cos \theta (x_2-1))}$. Here $k$ be the wavenumber and 
$\theta\in [ 0,\pi/2)$ is the incident angle. The total field $u_{\varepsilon}$, which consists of the incident wave $u^{i}$ and the scattered field $u_{\varepsilon}^s$ in $\Omega^{+}$ and the scattered field $u_{\varepsilon}^s$ only in $\Omega^{-}$, satisfies the Helmholtz equation
\begin{equation}\label{eq-Helmholtz}
\Delta u_{\varepsilon} + k^2 u_{\varepsilon} = 0   \quad\quad  \mbox{in} \; \Omega_{\varepsilon},
\end{equation}
and the boundary conditition
\begin{equation}\label{eq-bnd_cnd}
\dfrac{\partial u_{\varepsilon}}{\partial \nu} = 0  \quad \mbox{on} \; \partial \Omega_{\varepsilon}.
\end{equation}
Let $\kappa= k \sin \theta$, we look for quasi-periodic solutions such that $u_{\varepsilon}(x_1,x_2)=e^{i\kappa x_1}\tilde u_{\varepsilon}(x_1,x_2)$,
where $\tilde u_{\varepsilon}$ is a periodic function with $\tilde u_{\varepsilon}(x_1+d,x_2)=\tilde u_{\varepsilon}(x_1,x_2)$,
or equivalently,
\begin{equation}\label{eq-quasi_periodic}
u_{\varepsilon}(x_1+d,x_2)=e^{i\kappa d} u_{\varepsilon}(x_1,x_2).
\end{equation}
Define
$$
\kappa_n=\kappa+\dfrac{2\pi n}{d} \quad \mbox{and} \quad
\zeta_n(k)= \sqrt{k^2-\kappa_n^2},
$$
where the function $f(z)=\sqrt{z}$ is understood as an analytic function defined in the domain $\mathbf{C}\backslash\{-it: t\geq 0\}$ by
$$
z^{\frac{1}{2}} = |z|^{\frac{1}{2}} e^{\frac{1}{2} i\arg z}  
$$
throughout the paper.  Then it can be shown that the outgoing scattered field adopts the following Rayleigh-Bloch expansion (cf. \cite{bao95, bonnet_starling94, shipman10})
\begin{equation}\label{eq-rad_cond}
 u_{\varepsilon}^s(x_1,x_2) =  \sum_{n=-\infty}^{\infty} u_n^{s,+} e^{i \kappa_n x_1 + i\zeta_n  x_2 } \quad \mbox{and} \quad
 u_{\varepsilon}^s(x_1,x_2) =  \sum_{n=-\infty}^{\infty} u_n^{s,-} e^{i \kappa_n x_1 - i\zeta_n  x_2 }
\end{equation}
in the domain $\Omega^{+}$ and $\Omega^{-}$ respectively for some coefficients $u_n^{s,\pm}$.  
The expansion \eqref{eq-rad_cond} is usually referred to as the outgoing radiation condition and is imposed for the scattered field
above and below the metallic slab. 

Due to the quasi-periodicity of the solution, one can restrict $\kappa$ to the first Brillouin zone  $(-\pi/d,\pi/d]$.
Such $\kappa$ is called the reduced wave vector component \cite{bonnet_starling94, shipman10}.  
For given $k, \kappa, d$, we denote three sets of indices: $I_1=I_1(k, \kappa, d)= \{n: \abs{\kappa+2\pi n/d} < k\}$,  $I_2=I_2(k, \kappa, d)= \{n: 
\abs{\kappa+ 2\pi n/d} > k\}$, and $I_3=I_3(k, \kappa, d)= \{n: \abs{\kappa+2\pi n/d} = k\}$. The spacial harmonics $e^{i \kappa_n x_1 + i\zeta_n  x_2 }$ are called propagating modes, evanescent modes and linear modes for $n$ belonging to $I_1$, $I_2$ and $I_3$, respectively.

In contrast to the scattering problem for a single slit where one has
uniqueness, the solution to the scattering problem \eqref{eq-Helmholtz} - \eqref{eq-rad_cond} may not be unique. 
Indeed, the corresponding homogeneous problem with $u^i=0$ may attain non-trivial solutions $(k, u_\varepsilon)$ for $k\in\mathbf{R}$  
\cite{bonnet_starling94, shipman03, shipman07}.
Such real-valued $k$ is called singular frequency and $u_\varepsilon$ is the corresponding guided mode or surface bound state 
that decays exponentially away from the periodic structure \cite{bonnet_starling94, shipman03, shipman07}. 
On the other hand, there may also exist complex-valued $k$'s such that the homogeneous problem attains nontrivial solutions.
Such $k$ is called a resonance (or scattering resonance)  of the scattering problem, and the associated nontrivial solutions are called leaky modes (or quasi-normal modes). If the frequency of the incident wave is close the real part of the resonance (resonance frequency), then an enhancement of scattering is expected if the imaginary part of the resonance is small. This is the mechanism of resonant scattering. We refer to \cite{hai1, hai2, hai3, hai4, hai5, hai6, hai7} for the recent mathematical investigation of the interesting applications of resonances in super-resolution/super-focusing and metasurfaces. 
Finally, if $I_3(k, \kappa) \neq \emptyset$, then near a cut-off frequency where $k= \pm (\kappa + 2\pi n/d)$ for some integer $n$, 
the propagating spatial harmonics $e^{i \kappa_n x_1 \pm i\zeta_n  x_2 }$ becomes an evanescent mode or vice versa.
As such the diffracted field will exhibit anomalous behaviors near the cut-off frequencies and this is the so-called Rayleigh anomaly \cite{maystre12}. 
If resonances occur near such cut-off frequencies, then the enhancement behavior at resonant frequencies will be different from 
that of resonant phenomenon away from the cut-off frequencies.

In this paper, by using layer potential techniques, variational approaches and asymptotic analysis, we explore the field enhancement mechanism
 for the scattering problem in the diffraction regime. 
We derive the asymptotic expansions for both real eigenvalues and complex scattering resonances,
which lie below and above the light line $k=|\kappa|$ respectively.
It is known that surface bound states, for which the relation $k < |\kappa|$ holds, can not be excited by a plane incident wave $u^{i}= e^{i (\kappa x_1 - \zeta (x_2-1))}$ that satisfies $\kappa=k\sin\theta$. Whereas the quasi-modes, for which $|\kappa|<k$ holds, would be amplified if the incident frequency coincides the resonant frequencies. We analyze quantitatively the field enhancement at resonant frequencies, and show that
enhancement with an order of $O(\varepsilon^{-1})$ is achieved.
We also give the explicit dependence of the enhancement magnitude on the size of the period $d$ and the incident angle $\theta$.
Furthermore, the enhanced wave modes are also characterized in both near-field and far-field zones.
In this paper, we also demonstrate that if resonances are near the Rayleigh cut-off frequencies with a distance of $O(\varepsilon^{2\tau})$, where $0<\tau<1$, then the field enhancement at the resonant frequency becomes weaker and is of order $O(\varepsilon^{\tau-1})$.

The rest of the paper is organized as follows. We introduce layer potentials for the scattering problem and the boundary integral formulations 
in Section \ref{sec-bie}. The asymptotic expansion for the solution to the scattering problem is derived in Section \ref{sec-asy_sol}, which lays the foundation
for the quantitative analysis of field enhancement. The asymptotic expansions of real eigenvalues and complex resonances, and the investigation of field enhancement at resonant frequencies are presented in Section \ref{sec-eig_res}. In Section \ref{sec-res_Rayleigh}, we study the field enhancement 
when the resonant frequency is close to the Rayleigh cut-off frequencies. Finally, we discuss briefly on the embedded eigenvalues
associated with the scattering problem.

\setcounter{equation}{0}
\section{Boundary integral formulation}\label{sec-bie}
For each fixed $\kappa\in(-\pi/d,\pi/d]$, let $g^d(x, y)= g^d(x,y;\kappa)$ be the periodic Green's function solving the following equation:
$$ \Delta g^d(x,y) + k^2 g^d(x,y)=e^{i\kappa(x_1-y_1)} \sum_{n=-\infty}^{\infty} \delta(x_1-y_1-nd)\delta(x_2-y_2) \quad x, y\in\mathbf{R}^2.  $$
Then (cf. \cite{linton98})
\begin{equation}\label{eq-gd}
g^d(x,y; \kappa) = -\dfrac{i}{2d}  \sum_{n=-\infty}^{\infty} \dfrac{1}{\zeta_n(k)   } e^{i \kappa_n(x_1-y_1)+i\zeta_n(k)   |x_2-y_2| },
\end{equation}
where 
\begin{equation*}
\kappa_n=\kappa+\dfrac{2\pi n}{d} \quad \mbox{and} \quad
\zeta_n(k)   = \left\{
\begin{array}{lll}
\vspace*{5pt}
\sqrt{k^2-\kappa_n^2},  & \abs{\kappa_n} < k, \\
i\sqrt{\kappa_n^2-k^2},  & \abs{\kappa_n} > k. \\
\end{array}
\right.
\end{equation*}
We define the exterior Green's function with the Neumann boundary condition as $g^e(x,y)=g^e(x,y; \kappa) =g^d(x,y; \kappa)+g^d(x',y, \kappa)$, where 
\begin{equation*}
x' = \left\{
\begin{array}{ll}
(x_1, 2-x_2) & \mbox{if} \; x, y \in \Omega^+,  \\
(x_1,-x_2) &  \mbox{if} \; x, y \in \Omega^-.
\end{array}
\right.
\end{equation*}
It is easy to verify that $\dfrac{\partial g^e(x,y; \kappa)}{\partial \nu_y}=0$ on $\{y_2=1\}$ and $\{y_2=0\}$.

We also define the Green's function $g_\varepsilon^i(x,y)$ in the slit $S_{\varepsilon}^{(0)}$ with the Neumann boundary condition  as
$$ g_\varepsilon^i(x,y)= \sum_{m,n=0}^\infty c_{mn}\phi_{mn}(x)\phi_{mn}(y). $$
Here $c_{mn}=\dfrac{1}{k^2-(m\pi/\varepsilon)^2 - (n\pi)^2}$,  $\phi_{mn}=\sqrt{\dfrac{a_{mn}}{\varepsilon}}\cos\left(\dfrac{m\pi x_1}{\varepsilon}\right) \cos(n\pi x_2)$
and the coefficient $a_{mn}$ is given by
\begin{equation*}
a_{mn} = \left\{
\begin{array}{llll}
1  & m=n=0, \\
2  & m=0, n\ge 1 \quad \mbox{or} \quad n=0, m\ge 1, \\
4  & m\ge 1, n \ge 1.
\end{array}
\right.
\end{equation*}
It is easy to check that
$$ \Delta g_\varepsilon^i(x,y, k) + k^2 g_\varepsilon^i(x,y, k) = \delta(x-y), \quad x, y\in  S_{\varepsilon}^{(0)}. $$

\begin{figure}[!htbp]
\begin{center}
\includegraphics[height=5.6cm,width=14cm]{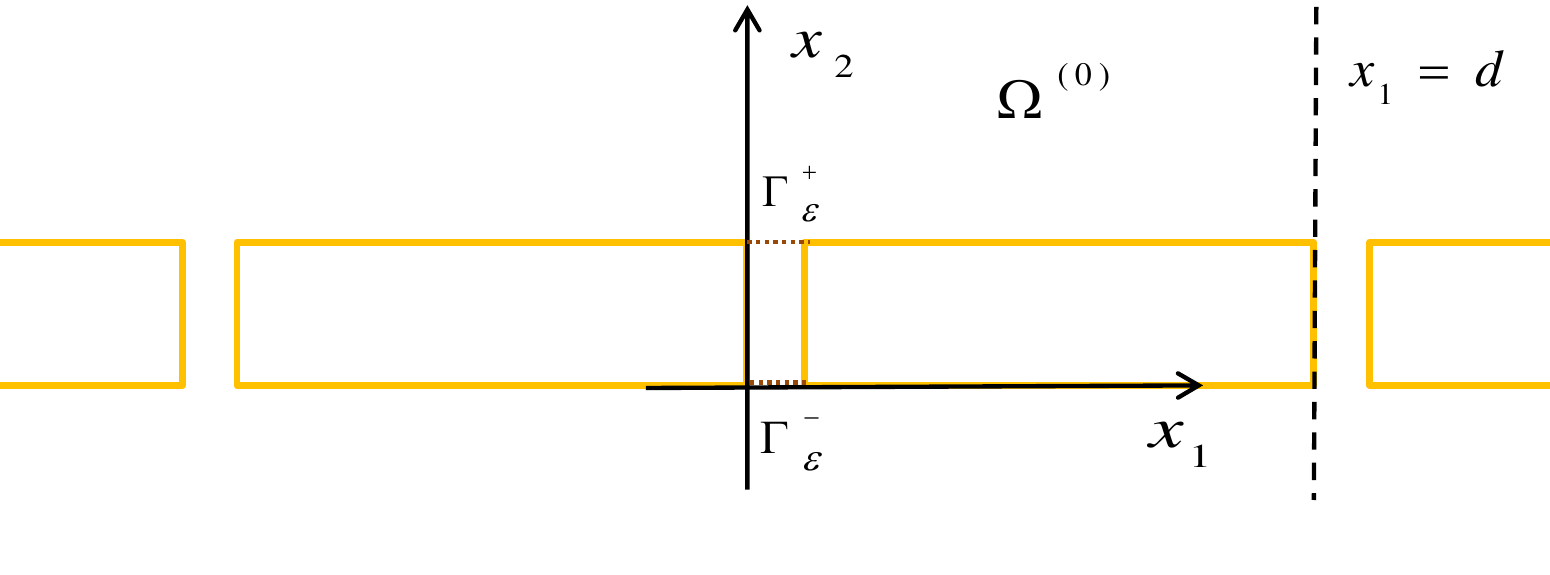}
\vspace{-10pt}
\caption{Problem geometry in one reference period $\Omega^{(0)}$.}\label{fig-prob_one_period}
\vspace{-10pt}
\end{center}
\end{figure}

To formulate the scattering problem as boundary integral equations, let us reduce the problem to one reference period $\Omega^{(0)}:=\{ x \in \mathbf{R}^2 \; | \;   0<x_1<d \}$ 
as shown in Figure \ref{fig-prob_one_period}.
Recall that the slit sitting in $\Omega^{(0)}$ is given by $S_\varepsilon^{(0)}:=\{(x_1,x_2)\;|\; 0<x_1<\varepsilon, 0<x_2<\ell \}$. We also denote the 
the upper and lower aperture of the slit by $\Gamma^{+}_\varepsilon$ and $\Gamma^{-}_\varepsilon$ respectively (see Figure \ref{fig-prob_one_period}).
\begin{lem}\label{eq-u_eps_formula}
Let $u_\varepsilon(x)$ be the solution of the scattering problem \eqref{eq-Helmholtz} - \eqref{eq-rad_cond}, then
\begin{eqnarray*}
u_\varepsilon(x) &=& \int_{\Gamma^+_\varepsilon} g^e(x,y) \dfrac{\partial u_\varepsilon(y)}{\partial y_2} ds_y + u^i+ u^r  \quad \mbox{for} \; x\in\Omega^{(0)} \cap \Omega^+. \\
u_\varepsilon(x) &=& -\int_{\Gamma^-_\varepsilon} g^e(x,y) \dfrac{\partial u_\varepsilon(y)}{\partial y_2} ds_y  \quad \mbox{for} \; x\in\Omega^{(0)} \cap \Omega^-. \\
u_\varepsilon(x) &=& \int_{\Gamma^-_\varepsilon} g_\varepsilon^i(x,y) \dfrac{\partial u_\varepsilon(y)}{\partial y_2} ds_y-\int_{\Gamma^+_\varepsilon} g_\varepsilon^i(x,y) \dfrac{\partial u_\varepsilon(y)}{\partial y_2} ds_y \quad \mbox{for} \; x\in S_\varepsilon^{(0)}.
\end{eqnarray*}
Here $u^r=e^{i (\kappa x_1 + \zeta (x_2-1))}$ is the reflected field of the ground plane $\{x_2=1\}$ without slits.
\end{lem}

\noindent\textbf{Proof}  For $x, y \in\Omega^{(0)}\cap \Omega^+$, define $\tilde g^e(x,y)=e^{-i\kappa(x_1-y_1)} g^e(x,y)$, and  $\tilde u^s_\varepsilon(y)=e^{-i\kappa y_1} u^s_\varepsilon(y)$. 
Then $\tilde g^e(x,y)$ and $\tilde u^s_\varepsilon(y)$ are periodic with respect to $y_1$.
Furthermore, a direction calculation yields
\begin{eqnarray*}
 \Delta \tilde g^e(x,y) - i2\kappa \dfrac{\partial  \tilde g^e(x,y)}{\partial y_1} + (k^2-\kappa^2) \tilde g^e(x,y)&=& \sum_{n=-\infty}^{\infty} \delta(x_1-y_1-nd)\delta(x_2-y_2), \\
 \Delta \tilde u_\varepsilon^s(y) + i2\kappa \dfrac{\partial  \tilde u_\varepsilon^s(y)}{\partial y_1} + (k^2-\kappa^2) \tilde u_\varepsilon^s(y)&=&0.  
\end{eqnarray*}
Choose $R>x$. Define the bounded domain above the perfect conductor by
$$ \Omega^+_R  := \{ x \in\Omega^{(0)}\; | \;   1<x_2 < R \}. $$
Denote the boundary of $\Omega^+_R $ by $\partial\Omega_R^{+}$. Let  $\Gamma_R^+:=\{ (x_1,R) \;|\;  0<x_1 < d \}$ and $\Gamma_1^+:=\{ (x_1,1) \;|\;  0<x_1 < d \}$ be its upper and lower boundaries
respectively.
An application of the Green's second identity and the divergence theorem over the domain $ \Omega^+_R$ yields
\begin{eqnarray*}
\tilde u_\varepsilon^s(x) &=& \int_{\Omega_R^{+}}  \Delta \tilde g^e(x,y) \tilde u_\varepsilon^s(y) -  \Delta  \tilde u_\varepsilon^s(y) \tilde g^e(x,y)
 -i2\kappa \int_{\Omega^+_R} \dfrac{\partial }{\partial y_1} \left(\tilde g^e(x,y) \tilde u_\varepsilon^s(y)\right) dy \\
 &=& \int_{\partial\Omega_R^{+}}   \dfrac{\partial \tilde g^e(x,y)}{\partial \nu_y}   \tilde u_\varepsilon^s(y)-\tilde g^e(x,y)  \dfrac{\partial  \tilde u_\varepsilon^s(y)}{\partial \nu_y}ds_y -i2\kappa \int_{\partial\Omega_R^{+}}(\tilde g^e(x,y) \tilde u_\varepsilon^s(y)) \nu_1 ds_y. 
\end{eqnarray*}
Using the fact that $\tilde g^e(x,y)$ and $\tilde u_\varepsilon^s(y)$ are periodic along $y_1$,
the second integral vanishes by noting that $\nu_1=0$ on the horizontal boundaries and $\nu$ take opposite signs on two vertical boundaries.
Substituting the Rayleigh expansions \eqref{eq-rad_cond} and \eqref{eq-gd}, it follows that
$$ \int_{\Gamma_R^+}  \dfrac{\partial \tilde g^e(x,y)}{\partial \nu_y}   \tilde u_\varepsilon^s(y)-\tilde g^e(x,y)  \dfrac{\partial  \tilde u_\varepsilon^s(y)}{\partial \nu_y}ds_y =0. $$
Thus, again applying the periodic boundary conditions for $\tilde g^e(x,y)$ and $\tilde u_\varepsilon^s(y)$, we have
$$ \tilde u_\varepsilon^s(x) =\int_{\Gamma^+_1}  \tilde g^e(x,y)  \dfrac{\partial  \tilde u_\varepsilon^s(y)}{\partial y_2} ds_y, $$
or equivalently,
\begin{eqnarray*}
\tilde u_\varepsilon (x) &=& \int_{\Gamma^+_1}  \tilde g^e(x,y)  \left(\dfrac{\partial \tilde u_\varepsilon(y)}{\partial y_2} -\dfrac{\partial  \tilde u^i(y)}{\partial y_2} \right) ds_y  + \tilde u^i (x), \\ 
&=& \int_{\Gamma^+_\varepsilon}  \tilde g^e(x,y)  \dfrac{\partial  \tilde u_\varepsilon(y)}{\partial y_2} ds_y -\int_{\Gamma^+_1}  \tilde g^e(x,y)  \dfrac{\partial  \tilde u^i(y)}{\partial y_2} ds_y  + \tilde u^i (x), \\ 
\end{eqnarray*}
where $\tilde u^i (x)=e^{-i\kappa x_1}u^i (x)$. A direct calculation leads to
$$ \int_{\Gamma^+_1}  \tilde g^e(x,y)  \dfrac{\partial  \tilde u^i(y)}{\partial y_2} ds_y = -e^{i\zeta(x_2-1)}. $$
Since $g^e(x,y)=e^{i\kappa(x_1-y_1)}\tilde g^e(x,y)$,  $u^s_\varepsilon(y)=e^{i\kappa y_1} \tilde u^s_\varepsilon(y)$, and $ u^i (x)=e^{i\kappa x_1}\tilde u^i (x)$.  It follows that
$$ u_\varepsilon(x) =\int_{\Gamma^+_\varepsilon}  g^e(x,y)  \dfrac{\partial  u_\varepsilon(y)}{\partial y_2} ds_y + u^i(x)+u^r(x), $$
where we have used the boundary condition $\partial_\nu u_\varepsilon=0$ on $\partial\Omega_\varepsilon$.

Following the similar lines above, it can be shown that 
$$ u_\varepsilon(x) = -\int_{\Gamma^-_\varepsilon} g^e(x,y) \dfrac{\partial u_\varepsilon}{\partial y_2} ds_y  \quad \mbox{for} \; x\in\Omega^{(0)} \cap \Omega^-. $$
On the other hand,  an application of Green's second identity over the slit  $S_\varepsilon^{(0)}$ gives rise to
$$  u_\varepsilon(x) = -\int_{\Gamma^+_\varepsilon\cup\Gamma^-_\varepsilon} g_\varepsilon^i(x,y) \dfrac{\partial u_\varepsilon}{\partial \nu_y} ds_y \quad \mbox{for} \; x\in S_\varepsilon^{(0)}. $$

\bigskip

Based upon Lemma  \ref{eq-u_eps_formula} and the continuity of the single layer potential (cf. \cite{kress}), we obtain the following boundary integral equations
defined over the slit apertures $\Gamma^\pm_\varepsilon$.
\begin{lem}
The following hold for the solution of the scattering problem \eqref{eq-Helmholtz} - \eqref{eq-rad_cond}:
\begin{equation}\label{eq-T1_e_eps}
u_\varepsilon(x) = \int_{\Gamma^+_\varepsilon}  g^e(x,y)\dfrac{\partial u_\varepsilon(y)}{\partial y_2} ds_y + u^i+u^r \quad \mbox{for} \; x\in\Gamma^+_\varepsilon.
\end{equation}
\begin{equation}\label{eq-T2_e_eps}
u_\varepsilon(x) = -\int_{\Gamma^-_\varepsilon}  g^e(x,y) \dfrac{\partial u_\varepsilon(y)}{\partial y_2} ds_y \quad \mbox{for} \; x\in\Gamma^-_\varepsilon.
\end{equation}
\begin{equation}\label{eq-T_i_eps}
u_\varepsilon(x) =  \int_{\Gamma^-_\varepsilon} g_\varepsilon^i(x,y) \dfrac{\partial u_\varepsilon(y)}{\partial y_2} ds_y-\int_{\Gamma^+_\varepsilon} g_\varepsilon^i(x,y) \dfrac{\partial u_\varepsilon(y)}{\partial y_2} ds_y
 \quad \mbox{for} \; x\in \Gamma^+_\varepsilon\cup\Gamma^-_\varepsilon.
\end{equation}
\end{lem}
An application of the above Lemma leads to the following system of integral equations:
\begin{equation} \label{eq-scattering2}
\left\{
\begin{array}{lll}
&\displaystyle{\int_{\Gamma^+_\varepsilon}  g^e(x,y)\dfrac{\partial u_\varepsilon(y)}{\partial y_2} ds_y 
+\int_{\Gamma^+_\varepsilon} g_\varepsilon^i(x,y) \dfrac{\partial u_\varepsilon(y)}{\partial y_2} ds_y
-\int_{\Gamma^-_\varepsilon} g_\varepsilon^i(x,y) \dfrac{\partial u_\varepsilon(y)}{\partial y_2} ds_y
+ u^i+u^r =0}, \quad \mbox{on} \,\, \Gamma^+_\varepsilon, \\ \\
& -\displaystyle{\int_{\Gamma^-_\varepsilon}  g^e(x,y)\dfrac{\partial u_\varepsilon(y)}{\partial y_2} ds_y 
+\int_{\Gamma^+_\varepsilon} g_\varepsilon^i(x,y) \dfrac{\partial u_\varepsilon(y)}{\partial y_2} ds_y
-\int_{\Gamma^-_\varepsilon} g_\varepsilon^i(x,y) \dfrac{\partial u_\varepsilon(y)}{\partial y_2} ds_y=0,} \quad \mbox{on} \,\, \Gamma^-_\varepsilon. 
\end{array}
\right.
\end{equation}

\begin{prop}
The scattering problem \eqref{eq-Helmholtz} - \eqref{eq-rad_cond} is equivalent to the system of boundary integral equations (\ref{eq-scattering2}). 
\end{prop}

It is clear that 
$$
\left.\dfrac{\partial u_\varepsilon}{\partial \nu}\right|_{\Gamma^+_\varepsilon} =\dfrac{\partial u_\varepsilon}{\partial y_2}(y_1, 1), 
\quad \left.\dfrac{\partial u_\varepsilon}{\partial \nu}\right|_{\Gamma^-_\varepsilon} =-\dfrac{\partial u_\varepsilon}{\partial y_2}(y_1, 0),\quad
(u^i+u^r)|_{\Gamma^+_\varepsilon}= 2e^{i\kappa x_1}.
$$
Note that the above functions are defined over narrow intervals with size $\varepsilon \ll 1$. To facilitate the analysis, we shall rescale the functions by introducing
$X=x_1/\varepsilon $ and $Y=y_1/\varepsilon$. Let us
define the following quantities: 
\begin{eqnarray*}
&& \varphi_1(Y):=  -\dfrac{\partial u_\varepsilon}{\partial y_2}( \varepsilon Y, 1); \\
&& \varphi_2(Y):=  \dfrac{\partial u_\varepsilon}{\partial y_2}(\varepsilon Y, 0);\\ 
&& f(X):= (u^i+u^r)(\varepsilon X, 1) = 2e^{i \kappa \varepsilon X};\\
&& G_\varepsilon^e(X, Y)=G_\varepsilon^e(X, Y, \kappa) :=  g^e(\varepsilon X,1; \varepsilon Y,1)=  g^e(\varepsilon X,0; \varepsilon Y,0)=-\dfrac{i}{d} \sum_{n=-\infty}^{\infty} \dfrac{1}{\zeta_n(k)   } e^{i \kappa_n\varepsilon(X-Y)}; \\
&& G_\varepsilon^i(X, Y) := g_\varepsilon^i(\varepsilon X, 1; \varepsilon Y, 1 ) = g_\varepsilon^i(\varepsilon X, 0; \varepsilon Y, 0 ) = \sum_{m,n=0}^\infty\dfrac{c_{mn}a_{mn}}{\varepsilon}\cos(m\pi X) \cos(m\pi Y); \\
&& \tilde G_\varepsilon^i(X, Y) := g_\varepsilon^i(\varepsilon X, 1; \varepsilon  Y, 0 ) = g_\varepsilon^i(\varepsilon  X, 0; \varepsilon  Y, 1 ) =\sum_{m,n=0}^\infty\dfrac{ (-1)^n c_{mn} a_{mn}}{\varepsilon}\cos(m\pi X) \cos(m\pi Y); 
\end{eqnarray*}
We also define three boundary integral operators:
\begin{eqnarray}\label{op_T}
  && (T^e \varphi)(X) = \int_0^1 G_\varepsilon^e(X, Y)  \varphi(Y) dY  \quad  X\in (0,1); \label{op_Te} \\
 && (T^i  \varphi) (X) = \int_0^1  G_\varepsilon^i(X, Y) \varphi(Y) dY  \quad  X\in (0,1); \label{op_Ti} \\
 && (\tilde T^i  \varphi)(X)  = \int_0^1 \tilde G_\varepsilon^i(X, Y) \varphi(Y) dY  \quad  X\in (0,1).  \label{op_tilde_Ti}
\end{eqnarray}
By a change of variable $x_1=\varepsilon X$ and $y_1=\varepsilon Y$ in \eqref{eq-scattering2}, the following proposition follows.
\begin{prop}
The system of equations (\ref{eq-scattering2}) is equivalent to the following one:
\begin{equation}\label{eq-scattering3}
\left[
\begin{array}{llll}
T^e+T^i    & \tilde{T}^i \\
\tilde{T}^i  & T^e+T^i 
\end{array}
\right] \left[
\begin{array}{llll}
\varphi_1    \\
\varphi_2 
\end{array}
\right]=\left[
\begin{array}{llll}
f/\varepsilon   \\
0
\end{array}
\right].
\end{equation}
\end{prop}

\section{Solution to the scattering problem}\label{sec-asy_sol}
In this section, based upon the integral equation formulation,
we derive the asymptotic expansion of the solution to the scattering problem \eqref{eq-Helmholtz} - \eqref{eq-rad_cond}.
From the derived expansion formulas, we obtain and classify different conditions for extraordinary field enhancement
for the underlying periodic structure. In addition, the expansions also lead to explicit asymptotic formulas for the enhanced wave fields 
presented in Section 4 and 5.

We begin with asymptotic expansions of the boundary integral operators $T^e$, $T^i$, and $\tilde{T}^i $. Then the solution
for the integral equation system \eqref{eq-scattering3} will be obtained, followed by the asymptotic expansion of the waves
in the far-field and near-field zones.

\subsection{Preliminaries}
We introduce several function spaces to be used in the rest of the paper.
Let $s \in \mathbf{R}$, we denote by $H^s(\mathbf{R})$ the standard fractional Sobolev space with the norm
$$
\|u\|_{H^s(\mathbf{R})}^2 = \int_{\mathbf{R}} (1+ |\xi|^2)^{s} |\hat{u}(\xi)|^2 d \xi,
$$
where $\hat{u}$ is the Fourier transform of $u$.  Let $I$ be a bounded open interval in  $\mathbf{R}$ and define
$$
H^s(I) := \{ u=U|_{I}  \;\big{|}\;  U \in H^s(\mathbf{R}) \}. 
$$
Then $H^s(I)$ is a Hilbert space with the norm
$$
\|u\|_{H^s(I)} = \inf \{ \|U\|_{H^s(\mathbf{R})} \;\big{|}\;  U \in H^s(\mathbf{R})\,\, \mbox {and}\,\, U|_I =u\}.
$$
We also define 
$$
\tilde H^s(I) := \{ u=U|_{I}  \;\big{|}\;  U \in H^s(\mathbf{R})\,\,  \mbox{and} \,\, supp \,U \subset \bar I   \}.
$$
One can show that the space $\tilde H^s(I)$ is the dual of $H^{-s}(I)$ and
the norm for $\tilde H^s(I)$ can be defined via the duality  \cite{adams}. As such $\tilde H^s(I)$  is also a Hilbert space. 
Here and henceforth, for simplicity, we denote $V_1 = \tilde H^{-\frac{1}{2}}(0, 1)$ and 
$V_2 =  H^{\frac{1}{2}}(0, 1)$. The duality between $V_1$ and $V_2$ will be denoted by 
$\langle u, v \rangle$ for any $u \in V_1$, $v \in V_2$. 

\subsection{Asymptotic expansion of the boundary integral operators}
For each fixed $\kappa\in(-\pi/d,\pi/d]$, we derive the asymptotic expansion of the boundary integral operators $T^e$, $T^i$, and $\tilde{T}^i $.
This is performed for $k$ away from the Rayleigh cut-off frequencies where $k=\kappa+2\pi n/d$ for some integer $n$. In this way we exclude the scenario where the Green's function $G_\varepsilon^e(X,Y)$ is not well defined because of a vanishing $\zeta_n(k)$. 
To this end, let us introduce 
a parameter $\tau$ such that $0\le\tau<1$
and denote
\begin{equation}\label{eq-delta}
\delta=O(\varepsilon^{2\tau}). 
\end{equation}
Let $B_{\delta}(z)$ be the disk with radius $\delta$ centered at $z$. 
Here and henceforth,
we define $$B_{\kappa,\delta}:=\bigcup_{n=-\infty}^{\infty} B_{\delta}(\kappa+2\pi n/d),$$ 
and restrict $k$ to the domain $\mathbf{R}^+\backslash B_{\kappa,\delta}$.
The following notations will be used throughout the paper.
\begin{eqnarray}
&& \beta_e(k,\kappa, d, \varepsilon)= \dfrac{1}{\pi} \left(\ln \varepsilon + \ln 2 + \ln\dfrac{\pi}{d}\right) + \left(\dfrac{1}{2\pi}\sum_{n\neq 0} \dfrac{1}{|n|} - \dfrac{i}{d} \sum_{n=-\infty}^{\infty}  \dfrac{1}{\zeta_n(k)}\right), \label{beta_e} \\
&& \beta_i(k, \varepsilon)= \dfrac{\cot k }{k \varepsilon} +  \dfrac{2\ln 2}{\pi},\label{beta_i} \\
&& \beta(k,\kappa, d, \varepsilon)=\beta_e(k,\kappa, d, \varepsilon) + \beta_i(k, \varepsilon), \label{beta} \\
&& \tilde \beta(k,\varepsilon) = \dfrac{1}{(k\sin k) \varepsilon}; \label{t_beta} \\
&& \rho(X, Y)= \dfrac{1}{\pi} \left[ \ln |X-Y| + \ln \left(\abs{\sin \left(\frac{\pi(X-Y)}{2}\right)}\right)+ \ln \left(\abs{\sin \left(\frac{\pi(X+Y)}{2}\right)}\right)\right]\label{kernel_k}
\end{eqnarray}

\medskip
\noindent\textbf{Remark 3.1}
\textit{From the definition of $\zeta_n (k)$, it is clear that  $\zeta_n(k)  \neq 0$ for $k\in\mathbf{R}^+\backslash B_{\kappa,\delta}$. 
In the above and throughout, the substraction
 $$\displaystyle{\frac{1}{2\pi}\sum_{n\neq 0} \frac{1}{|n|} - \frac{i}{d} \sum_{n=-\infty}^{\infty}\frac{1}{\zeta_n(k)  }}$$  is viewed as the sum of the convergent series 
 $$\displaystyle{\sum_{n\neq 0} \left(\frac{1}{2\pi}\frac{1}{|n|} - \frac{i}{d}\frac{1}{\zeta_n(k)}\right) - \frac{i}{d}\frac{1}{\zeta_0(k) }}. $$ Therefore, the scalar function $\beta_e(k,\kappa, d, \varepsilon)$ is well defined.} \\

\medskip
\noindent\textbf{Remark 3.2}
\textit{The real number $\delta=O(\varepsilon^{2\tau})$ denotes the distance from the Rayleigh cut-off frequencies.
When $\tau=0$ and $\delta=O(1)$, the frequency in $\mathbf{R}^+\backslash B_{\kappa,\delta}$ is then far away from the Rayleigh cut-off frequencies. 
On the other hand, if $0<\tau<1$, a frequency in $\mathbf{R}^+\backslash B_{\kappa,\delta}$ can be close to the Rayleigh anomaly. 
It should be pointed out that the assumption
$0\le\tau<1$ is essential for a uniform asymptotic expansions of the Green's functions, which are given in Lemma \ref{lem-green} below.
The treatment for $\tau>1$ would require more dedicated asymptotic analysis. A promising approach is to use new Green's functions that converge
rapidly around the Rayleigh cut-off frequencies \cite{bruno14}.} \\

The asymptotic expansions for the kernels $G_\varepsilon^e$,  $G_\varepsilon^i$, 
and $\tilde G_\varepsilon^i $ are given in the following Lemma.
\begin{lem} \label{lem-green} 
For $\kappa\in(-\pi/d,\pi/d]$ and $k\in\mathbf{R}^+\backslash B_{\kappa,\delta}$,  where $\delta=O(\varepsilon^{2\tau})$ and $0\le\tau<1$. If $k\varepsilon \ll 1$, then
\begin{eqnarray}\label{eq-G_asy}
G_\varepsilon^e(X, Y)&=&G_\varepsilon^e(X, Y; \kappa)  = \beta_e(k,\kappa, d, \varepsilon) + \dfrac{1}{\pi}  \ln |X-Y| +  r_1^{\varepsilon}(X, Y), \label{Ge_exp} \\
G_\varepsilon^i(X, Y) 
 &=& \beta_i(k ,\varepsilon) + \dfrac{1}{\pi} \left[ \ln \left(\abs{\sin \left(\frac{\pi(X+Y)}{2}\right)}\right) + \ln \left(\abs{\sin \left(\frac{\pi(X-Y)}{2}\right)}\right) \right] \nonumber\\
 && +  r_2^{\varepsilon}(X, Y), \label{Gi_exp} \\ 
\tilde G_\varepsilon^i(X, Y) &= &  \tilde \beta(k,\varepsilon) + \tilde \rho_\infty(X, Y). \label{tGi_exp}
\end{eqnarray}
Here $r_1^{\varepsilon}(X, Y)$, $r_2^{\varepsilon}(X, Y)$, and  $\tilde \rho_\infty(X, Y)$ are bounded functions with
$r_1^{\varepsilon}\sim O(k\varepsilon^{1-\tau})$, $r_2^{\varepsilon}\sim O((k\varepsilon)^2)$, and 
$\tilde \rho_\infty\sim O(\exp(-1/\varepsilon))$ for all $X, Y\in(0,1)$.
\end{lem}

\noindent\textbf{Proof} The asymptotic expansion for the kernels $G_\varepsilon^i$ and $\tilde G_\varepsilon^i $ has been derived in \cite{lin_zhang16}. See also Appendix A for the proof.
To obtain asymptotic expansion for the kernel $G_\varepsilon^e$, let us consider the case when $\kappa=0$.
The case of $\kappa\neq0$ can be calculated in a similar fashion.
Using Taylor expansion, we have
\begin{eqnarray*}
\sum_{|\kappa_n|>k} \dfrac{1}{\zeta_n(k)   } e^{i \kappa_n\varepsilon(X-Y)} &=& -\dfrac{id}{2\pi} \sum_{|\kappa_n|>k}\dfrac{1}{|n|\sqrt{1-(kd/2\pi n)^2}}e^{i \frac{2\pi n}{d}\varepsilon(X-Y)}  \\
 &=& -\dfrac{id}{2\pi} \sum_{|\kappa_n|>k}\dfrac{1}{|n|} \left( 1+\sum_{m=1}^\infty \dfrac{1\cdot3\cdots (2m-1)}{2^m m!} \left(\dfrac{kd}{2\pi n}\right)^{2m} \right) e^{i \frac{2\pi n}{d}\varepsilon(X-Y)}.
\end{eqnarray*}
Then it follows that
\begin{eqnarray}\label{eq-sum_fourier}
\sum_{n=-\infty}^{\infty} \dfrac{1}{\zeta_n(k)   } e^{i \kappa_n\varepsilon(X-Y)} &=&\sum_{|\kappa_n|<k} \dfrac{1}{\zeta_n(k)   } e^{i \kappa_n\varepsilon(X-Y)} + \sum_{|\kappa_n|>k} \dfrac{1}{\zeta_n(k)   } e^{i \kappa_n\varepsilon(X-Y)}  \nonumber \\
&=&  -\dfrac{id}{2\pi} \sum_{n\neq0}\dfrac{1}{|n|} e^{i \frac{2\pi n}{d}\varepsilon(X-Y)} 
+ \dfrac{id}{2\pi} \sum_{|\kappa_n|<k, n\neq0}\dfrac{1}{|n|}e^{i \frac{2\pi n}{d}\varepsilon(X-Y)} + \sum_{|\kappa_n|<k} \dfrac{1}{\zeta_n(k)   } e^{i \kappa_n\varepsilon(X-Y)}    \nonumber \\
&& -\dfrac{id}{2\pi}  \sum_{m=1}^\infty \dfrac{1\cdot3\cdots (2m-1)}{2^m m!} \sum_{|\kappa_n|>k} \left(\dfrac{kd}{2\pi n}\right)^{2m}\dfrac{1}{|n|} e^{i \frac{2\pi n}{d}\varepsilon(X-Y)}.
\end{eqnarray}
Notice that
\begin{eqnarray*}
&& -\sum_{n\neq0}\dfrac{1}{|n|} e^{i \frac{2\pi n}{d}\varepsilon(X-Y)}= \ln\left(4\sin^2\dfrac{\pi \varepsilon(X-Y)}{d}\right), \\
&& \sum_{|\kappa_n|<k, n\neq 0}\dfrac{1}{|n|}e^{i \frac{2\pi n}{d}\varepsilon(X-Y)} =  \sum_{|\kappa_n|<k, n\neq 0}\dfrac{1}{|n|} + O(k\varepsilon(X-Y)), \\
&& \sum_{|\kappa_n|<k}\dfrac{1}{\zeta_n(k) }e^{i \frac{2\pi n}{d}\varepsilon(X-Y)} =   \sum_{|\kappa_n|<k} \dfrac{1}{\zeta_n(k)   } + O(k\varepsilon^{1-\tau}(X-Y)).
\end{eqnarray*}
Furthermore, for $m\ge1$,
\begin{eqnarray*}
&& \sum_{|\kappa_n|>k}\dfrac{1}{|n|^{2m+1}}e^{i \frac{2\pi n}{d}\varepsilon(X-Y)} \\
&=& \sum_{n\neq 0 }\dfrac{1}{|n|^{2m+1}}e^{i \frac{2\pi n}{d}\varepsilon(X-Y)} -  \sum_{|\kappa_n|<k, n\neq 0}\dfrac{1}{|n|^{2m+1}}e^{i \frac{2\pi n}{d}\varepsilon(X-Y)} \\
&=& \sum_{n\neq 0 }\dfrac{1}{|n|^{2m+1}}+ O(\varepsilon^{2m}(X-Y)^{2m})\ln (\varepsilon(X-Y)) - \sum_{|\kappa_n|<k, n\neq0}\dfrac{1}{|n|^{2m+1}} + O(k\varepsilon(X-Y)) \\
&=& \sum_{|\kappa_n|>k}\dfrac{1}{|n|^{2m+1}} + O(k\varepsilon(X-Y)).
\end{eqnarray*}
Substituting the above into \eqref{eq-sum_fourier}, we have
\begin{eqnarray*}
\sum_{n=-\infty}^{\infty} \dfrac{1}{\zeta_n(k)   } e^{i \kappa_n\varepsilon(X-Y)} &=& \dfrac{id}{2\pi}\ln\left(4\sin^2\dfrac{\pi \varepsilon(X-Y)}{d}\right) + \dfrac{id}{2\pi} \sum_{|\kappa_n|<k, n\neq 0}\dfrac{1}{|n|} + \sum_{|\kappa_n|<k} \dfrac{1}{\zeta_n(k)   } \\
&& -\dfrac{id}{2\pi}  \sum_{m=1}^\infty \dfrac{1\cdot3\cdots (2m-1)}{2^m m!} \sum_{|\kappa_n|>k}  \left(\dfrac{kd}{2\pi}\right)^{2m}\dfrac{1}{|n|^{2m+1}} + O(k\varepsilon^{1-\tau}(X-Y)).
\end{eqnarray*}
Using the relation  
$$ \sum_{|\kappa_n|>k} \dfrac{1}{|n|\sqrt{1-(kd/2\pi n)^2}}-\dfrac{1}{|n|}=\sum_{|\kappa_n|>k} \sum_{m=1}^\infty \dfrac{1\cdot3\cdots (2m-1)}{2^m m!}  \left(\dfrac{kd}{2\pi}\right)^{2m}\dfrac{1}{|n|^{2m+1}} ,$$
we get
\begin{eqnarray*}
\sum_{n=-\infty}^{\infty} \dfrac{1}{\zeta_n(k)   } e^{i \kappa_n\varepsilon(X-Y)} &=& \dfrac{id}{2\pi}\ln\left(4\sin^2\dfrac{\pi \varepsilon(X-Y)}{d}\right) + \left(\dfrac{id}{2\pi}\sum_{n\neq 0} \dfrac{1}{|n|} + \sum_{n=-\infty}^{\infty}  \dfrac{1}{\zeta_n(k)   }\right)+ O(k\varepsilon^{1-\tau}(X-Y)).
\end{eqnarray*}
Therefore, the desired asymptotic expansion follows by noting that 
$$\displaystyle{G_\varepsilon^e(X,Y) = -\dfrac{i}{d}\sum_{n=-\infty}^{\infty} \dfrac{1}{\zeta_n(k)   } e^{i \kappa_n\varepsilon(X-Y)}}.$$ \qed

Let $\rho(X,Y)$ and $\tilde \rho_{\infty}(X,Y)$ be defined in \eqref{kernel_k} and \eqref{tGi_exp} respectively.
Set $\rho_{\infty}(X,Y)= r_1^{\varepsilon}(X,Y)+ r_2^{\varepsilon}(X,Y)$,  where $ r_1^{\varepsilon}(X,Y)$ and $ r_2^{\varepsilon}(X,Y)$ are given in
\eqref{Ge_exp} and \eqref{Gi_exp} respectively.
We denote by $K$, $K_{\infty}$, $\tilde K_{\infty}$ the integral operators corresponding to the Schwarz kernels $\rho(X,Y)$, $\rho_{\infty}(X,Y)$ and $\tilde \rho_{\infty}(X,Y)$, respectively. We also define the operator $P: V_1 \to V_2$ by 
$$
P \varphi(X) = \langle \varphi, 1 \rangle 1,
$$
where $1$ is a function defined on the interval $(0, 1)$ and is equal to one therein.  We will use this notation in the sequel. One can easily check that $1 \in V_2$. Thus the above definition is valid.

\begin{lem} \label{lem-operators}
For $\kappa\in(-\pi/d,\pi/d]$ and $k\in\mathbf{R}^+\backslash B_{\kappa,\delta}$, where $\delta=O(\varepsilon^{2\tau})$ and $0\le\tau<1$.
\begin{enumerate}
\item[(1)]
The operator $T^e+T^i $ admits the following decomposition:
$$ T^e+T^i = \beta P + K + K_\infty. $$
Moreover, $K_\infty$ is bounded from $V_1$ to $V_2$ with operator norm
$\| K_{\infty}\|  \lesssim \varepsilon^{1-\tau}$ uniformly for bounded $k$'s.

\item[(2)]
The operator  $\tilde T^i$ admits the following decomposition: 
$$\tilde T^i = \tilde \beta P + \tilde K_\infty,$$ 
Moreover, $\tilde K_\infty$ is bounded from $V_1$ to $V_2$ with operator norm
$\| K_{\infty}\|  \lesssim e^{-1/\varepsilon}$ uniformly for bounded $k$'s. 

\item[(3)]
The operator $K$ is bounded from $V_1$ to $V_2$ with a bounded inverse.
Moreover, 
$$\alpha:=\langle K^{-1} 1, 1 \rangle_{L^2(0,1)} \neq 0. $$
\end{enumerate}
\end{lem}

\noindent The proof of (1) and (2)  follows directly from the definition of the operators $T^e$, $T^i$ and $\tilde T^i $ in \eqref{op_Te} - \eqref{op_tilde_Ti}
and the asymptotic expansions of their kernels (see Lemma \ref{lem-green}).
The proof of (3) can be found in Theorem 4.1 and Lemma 4.2 of \cite{eric10}.

\subsection{ Asymptotic expansion of the solution to the system \eqref{eq-scattering3}}
Define
\begin{equation*}
\mathbb{P}= \left[
\begin{array}{cc}
\beta P    & \tilde \beta P \\
\tilde \beta P  & \beta P 
\end{array}
\right],
\quad
\mathbb{K}_\infty=
\left[
\begin{array}{cc}
K_\infty    & \tilde K_\infty \\
\tilde K_\infty  & K_\infty
\end{array}
\right], 
\quad
\mathbf{f}=\left[ \begin{array}{c}
f/\varepsilon \\
0
\end{array}\right],
\quad\mbox{and}\quad
\mathbb{L}=K \mathbb{I}  +\mathbb{K}_\infty.
\end{equation*}
Then from the decomposition of the operators in Lemma \ref{lem-operators}, we may rewrite the system of the integral equations
\eqref{eq-scattering3} as
\begin{equation}\label{eq-scattering4}
(\mathbb{P} + \mathbb{L}) \boldsymbol{\varphi} =   \mathbf{f}. 
\end{equation}
Next, we derive the asymptotic expansion of the solution $\boldsymbol{\varphi}$.
By Lemma \ref{lem-operators}, it is also easy to see that $\mathbb{L}$ is invertible for sufficiently small 
$\varepsilon$. 
Applying the Neumann series yields
$$ \mathbb{L}^{-1}  = \left(K\mathbb{I}+\mathbb{K}_\infty\right)^{-1} = \left(\sum_{j=0}^\infty (-1)^j\left(K^{-1}\mathbb{K_\infty} \right)^j \right)  K^{-1} = K^{-1}\mathbb{I} + O\left(k\varepsilon^{1-\tau}\right). $$
Therefore, the following lemma follows.
\begin{lem} \label{lem-L_inv} 
 Let $\mathbf{e}_1 = [1, 0]^T $ and $\mathbf{e}_2 = [0, 1]^T $. Then
\begin{equation*} \label{eq-II}
\mathbb{L}^{-1}  \mathbf{e}_1 = K^{-1}1 \cdot \mathbf{e}_1 + O(k\varepsilon^{1-\tau}), \quad
\mathbb{L}^{-1}  \mathbf{e}_2 = K^{-1}1 \cdot \mathbf{e}_2 + O(k\varepsilon^{1-\tau}),
\end{equation*}
and 
\begin{equation} \label{eq-III}
\langle \mathbb{L}^{-1}  \mathbf{e}_1,  \mathbf{e}_1 \rangle = \alpha +  O(k\varepsilon^{1-\tau}), \quad 
\langle \mathbb{L}^{-1}  \mathbf{e}_1,  \mathbf{e}_2\rangle =O(k\varepsilon^{1-\tau}).
\end{equation}
\end{lem}
The following identities are proved in \cite{lin_zhang16}. 
\begin{lem} \label{lem-identity1}
 Let $\mathbf{e}_1 = [1, 0]^T $ and $\mathbf{e}_2 = [0, 1]^T $. Then
\begin{equation*} 
\label{eq-l} \langle \mathbb{L}^{-1}  \mathbf{e}_1,  \mathbf{e}_1 \rangle = \langle \mathbb{L}^{-1}  \mathbf{e}_2,  \mathbf{e}_2 \rangle, \quad \langle \mathbb{L}^{-1}  \mathbf{e}_1,  \mathbf{e}_2 \rangle = \langle   \mathbb{L}^{-1}\mathbf{e}_2,  \mathbf{e}_1 \rangle.
\end{equation*}
\end{lem}

\medskip

By applying $\mathbb{L}^{-1}$ on both sides of \eqref{eq-scattering4},  we see that
\begin{equation}\label{eq-scattering5}
\mathbb{L}^{-1} \;\mathbb{P} \; \boldsymbol{\varphi} + \boldsymbol{\varphi} =  \mathbb{L}^{-1}  \mathbf{f}.
\end{equation}
Note that
 $$  \mathbb{P} \; \boldsymbol{\varphi} = \beta  \langle  \boldsymbol{\varphi}, \mathbf{e}_1 \rangle \mathbf{e}_1  + \beta \langle  \boldsymbol{\varphi}, \mathbf{e}_2 \rangle \mathbf{e}_2
 + \tilde \beta \langle  \boldsymbol{\varphi}, \mathbf{e}_2 \rangle \mathbf{e}_1  + \tilde \beta \langle  \boldsymbol{\varphi}, \mathbf{e}_1 \rangle \mathbf{e}_2,$$
the above operator equation can be written as
\begin{equation}\label{Op_eqns_1}
\beta \langle  \boldsymbol{\varphi}, \mathbf{e}_1 \rangle  \mathbb{L}^{-1}  \mathbf{e}_1  + \beta  \langle  \boldsymbol{\varphi}, \mathbf{e}_2 \rangle  \mathbb{L}^{-1}  \mathbf{e}_2 + 
\tilde \beta \langle  \boldsymbol{\varphi}, \mathbf{e}_2 \rangle \mathbb{L}^{-1} \mathbf{e}_1  + \tilde \beta \langle  \boldsymbol{\varphi}, \mathbf{e}_1 \rangle \mathbb{L}^{-1} \mathbf{e}_2+
\boldsymbol{\varphi} =   \mathbb{L}^{-1}  \mathbf{f}. 
\end{equation}
By taking the inner product of \eqref{Op_eqns_1} with $\mathbf{e}_1$ and $\mathbf{e}_2$ respectively, it follows that
\begin{equation}\label{eqn-linear_sys}
( \mathbb{M}+\mathbb{I})\left[
 \begin{array}{llll}
\langle  \boldsymbol{\varphi}, \mathbf{e}_1 \rangle  \\
\langle  \boldsymbol{\varphi}, \mathbf{e}_2 \rangle
\end{array}
\right] =
\left[
\begin{array}{llll}
 \langle \mathbb{L}^{-1} \mathbf{f}, \mathbf{e}_1 \rangle   \\
 \langle \mathbb{L}^{-1} \mathbf{f}, \mathbf{e}_2 \rangle 
\end{array}
\right],
\end{equation}
where the matrix $\mathbb{M}$ is defined as
\begin{equation}  \label{eqn-matrix-m}
\mathbb{M}:=
\beta  \left[
\begin{array}{llll}
 \langle \mathbb{L}^{-1}  \mathbf{e}_1,  \mathbf{e}_1 \rangle  &  \langle \mathbb{L}^{-1}  \mathbf{e}_2,  \mathbf{e}_1\rangle \\
 \langle \mathbb{L}^{-1}  \mathbf{e}_1,  \mathbf{e}_2 \rangle  &  \langle \mathbb{L}^{-1}  \mathbf{e}_2,  \mathbf{e}_2 \rangle 
 \end{array}
\right] 
+
\tilde \beta \left[
\begin{array}{llll}
 \langle \mathbb{L}^{-1}  \mathbf{e}_2,  \mathbf{e}_1 \rangle  &  \langle \mathbb{L}^{-1}  \mathbf{e}_1,  \mathbf{e}_1\rangle \\
 \langle \mathbb{L}^{-1}  \mathbf{e}_2,  \mathbf{e}_2 \rangle  &  \langle \mathbb{L}^{-1}  \mathbf{e}_1,  \mathbf{e}_2 \rangle 
 \end{array}
\right]. 
\end{equation}
Therefore,
\begin{equation}\label{eq-phi_dot_e1e2}
\left[
\begin{array}{llll}
\langle  \boldsymbol{\varphi}, \mathbf{e}_1 \rangle  \\
\langle  \boldsymbol{\varphi}, \mathbf{e}_2 \rangle
\end{array}
\right] =( \mathbb{M}+\mathbb{I})^{-1} 
\left[
\begin{array}{llll}
 \langle \mathbb{L}^{-1} \mathbf{f}, \mathbf{e}_1 \rangle   \\
 \langle \mathbb{L}^{-1} \mathbf{f}, \mathbf{e}_2 \rangle 
\end{array}
\right].
\end{equation}
Substituting into \eqref{eq-scattering5} yields 
\begin{equation}\label{eq-phi}
\boldsymbol{\varphi} = \mathbb{L}^{-1}  \mathbf{f}
- \bigg[  \mathbb{L}^{-1}  \mathbf{e}_1 \quad  \mathbb{L}^{-1}  \mathbf{e}_2 \bigg]  
\left[
\begin{array}{llll}
\beta  &   \tilde \beta \\
\tilde \beta  &\beta
\end{array}
\right] 
(\mathbb{M}+\mathbb{I})^{-1} 
\left[
\begin{array}{llll}
 \langle \mathbb{L}^{-1} \mathbf{f}, \mathbf{e}_1 \rangle   \\
 \langle \mathbb{L}^{-1} \mathbf{f}, \mathbf{e}_2 \rangle 
\end{array}
\right].
\end{equation}

From Lemma \ref{lem-identity1}, it is observed that
\begin{equation*}
\mathbb{M} = 
\left(\beta  +
\tilde \beta 
\left[
\begin{array}{llll}
0 & 1 \\
1 & 0 
\end{array}
\right]\right)
\left[
\begin{array}{llll}
 \langle \mathbb{L}^{-1}  \mathbf{e}_1,  \mathbf{e}_1 \rangle  &  \langle \mathbb{L}^{-1}  \mathbf{e}_1,  \mathbf{e}_2\rangle \\
 \langle \mathbb{L}^{-1}  \mathbf{e}_1,  \mathbf{e}_2 \rangle  &  \langle \mathbb{L}^{-1}  \mathbf{e}_1,  \mathbf{e}_1 \rangle 
 \end{array}
\right].
\end{equation*}
A straightforward calculation shows that the eigenvalues of $\mathbb{M}+\mathbb{I}$ are
\begin{eqnarray}\label{eq-eigen_M}
\lambda_1(k; \kappa, d, \varepsilon) &=& 1+(\beta+ \tilde \beta )   \left(\langle \mathbb{L}^{-1}  \mathbf{e}_1,  \mathbf{e}_1 \rangle  + \langle \mathbb{L}^{-1}  \mathbf{e}_1,  \mathbf{e}_2\rangle\right), \label
{eq-lambda1}\\
\lambda_2(k; \kappa, d, \varepsilon) &=& 1+(\beta-\tilde \beta )   \left(\langle \mathbb{L}^{-1}  \mathbf{e}_1,  \mathbf{e}_1 \rangle  - \langle \mathbb{L}^{-1}  \mathbf{e}_1,  \mathbf{e}_2\rangle\right),
\label{eq-lambda2}
\end{eqnarray}
and the associated eigenvectors are $[1 \quad 1]^T$ and $[1 \quad -1]^T$.  
For simplicity of notation, we define two scalar functions
\begin{equation}\label{eq-def_pq}
p(k; \kappa, d, \varepsilon):=\varepsilon \lambda_1(k; \kappa, d, \varepsilon) \quad \mbox{and} \quad q(k; \kappa, d, \varepsilon):=\varepsilon \lambda_2(k; \kappa, d, \varepsilon).
\end{equation}
We also define
\begin{equation}\label{eq-gamma}
\gamma(k, \kappa, d)=\dfrac{1}{\pi} \left(3\ln 2 + \ln\dfrac{\pi}{d}\right) + \left(\dfrac{1}{2\pi}\sum_{n\neq 0} \dfrac{1}{|n|} - \dfrac{i}{d} \sum_{n=-\infty}^{\infty}  \dfrac{1}{\zeta_n(k) } \right).
\end{equation}
Then a combination of  \eqref{eq-lambda1} - \eqref{eq-def_pq}, the expressions \eqref{beta_e} - \eqref{t_beta} for $\beta$ and $\tilde\beta$, 
and Lemma \ref{lem-L_inv} yields
\begin{equation}\label{eq-formula_p}
p(k; \kappa, d, \varepsilon)=\varepsilon +\left[ \dfrac{\cot k }{k} + \dfrac{1}{k\sin k} +  \varepsilon \gamma(k, \kappa, d) + \dfrac{1}{\pi} \varepsilon \ln \varepsilon   \right]   \left(\alpha + O(k\varepsilon^{1-\tau}) \right),
\end{equation}
and
\begin{equation}\label{eq-formula_q}
q(k; \kappa, d, \varepsilon)=\varepsilon +\left[ \dfrac{\cot k }{k} - \dfrac{1}{k\sin k} +  \varepsilon \gamma(k, \kappa, d) + \dfrac{1}{\pi} \varepsilon \ln \varepsilon   \right]   \left(\alpha +O(k\varepsilon^{1-\tau}) \right),
\end{equation}

\begin{lem} \label{lem-phi}
Let $\kappa= k\sin \, \theta$, $\delta=O(\varepsilon^{2\tau})$ where $0\le\tau<1$,  and 
$k\in\mathbf{R}^+\backslash B_{\kappa,\delta}$ be bounded and not an eigenvalue of the scattering operator.
Then the following asymptotic expansion holds for the solution $\boldsymbol{\varphi}$ of \eqref{eq-scattering3} in $V_1 \times V_1$:
\begin{eqnarray} \label{eq-varphi}
 \boldsymbol{\varphi} &=& K^{-1}1 
 \cdot  \left[ \kappa \cdot O(1) \cdot\mathbf{e}_1 + \dfrac{\alpha}{p}  (\mathbf{e}_1 + \mathbf{e}_2 ) + \dfrac{\alpha}{q} (\mathbf{e}_1 - \mathbf{e}_2)  \right] \nonumber + \left(\dfrac{\alpha}{p} +  \dfrac{\alpha}{q} \right) \cdot O(k\varepsilon^{1-\tau}) +  O(k\varepsilon^{1-\tau}).
\end{eqnarray}
Moreover, 
\begin{equation}  \label{eq-varphi-1}
\left[
\begin{array}{llll}
\langle  \boldsymbol{\varphi}, \mathbf{e}_1 \rangle  \\
\langle  \boldsymbol{\varphi}, \mathbf{e}_2 \rangle
\end{array}
\right]  
= \bigg[\alpha+O(\varepsilon^{1-\tau}) \bigg] \left(
\dfrac{1}{p}
\left[
\begin{array}{cc}
1  \\
1 
\end{array}
\right] 
+
\dfrac{1}{q}
\left[
\begin{array}{cc}
1 \\
-1 
\end{array}
\right] \right).
\end{equation}
\end{lem}

\noindent\textbf{Proof}.  The matrix $\mathbb{M}+\mathbb{I}$ has two eigenvalues 
$\lambda_1$ and $\lambda_2$ given by \eqref{eq-lambda1} and \eqref{eq-lambda2}, 
which are associated with the eigenvectors $[1 \quad 1]^T$ and $[1 \quad -1]^T$ respectively. 
Thus
\begin{equation*}
(\mathbb{M}+\mathbb{I})^{-1}  =
\dfrac{1}{2\lambda_1}
\left[
\begin{array}{cc}
1 & 1 \\
1 & 1 
\end{array}
\right] 
+
\dfrac{1}{2\lambda_2}
\left[
\begin{array}{cc}
1 & -1 \\
-1 & 1 
\end{array}
\right].
\end{equation*}
By substituting into \eqref{eq-phi} \eqref{eq-phi_dot_e1e2},  it follows that 
\begin{eqnarray*}
\left[
\begin{array}{llll}
\langle  \boldsymbol{\varphi}, \mathbf{e}_1 \rangle  \\
\langle  \boldsymbol{\varphi}, \mathbf{e}_2 \rangle
\end{array}
\right] 
&=&
\dfrac{1}{2\lambda_1}\langle \mathbb{L}^{-1} \mathbf{f}, \mathbf{e}_1+\mathbf{e}_2 \rangle
\left[
\begin{array}{cc}
1  \\
1 
\end{array}
\right] 
+
\dfrac{1}{2\lambda_2(k,\varepsilon)}\langle \mathbb{L}^{-1} \mathbf{f}, \mathbf{e}_1-\mathbf{e}_2 \rangle
\left[
\begin{array}{cc}
1 \\
-1 
\end{array}
\right],
\end{eqnarray*}
and 
\begin{eqnarray*}
\boldsymbol{\varphi} &=& \mathbb{L}^{-1}  \mathbf{f}
- \dfrac{1}{2\lambda_1}\bigg[  \mathbb{L}^{-1}  \mathbf{e}_1 \quad  \mathbb{L}^{-1}  \mathbf{e}_2 \bigg]  
\left[
\begin{array}{llll}
\beta  &   \tilde \beta \\
\tilde \beta &\beta 
\end{array}
\right] 
\left[
\begin{array}{cc}
1 & 1 \\
1 & 1 
\end{array}
\right] 
\left[
\begin{array}{llll}
 \langle \mathbb{L}^{-1} \mathbf{f}, \mathbf{e}_1 \rangle   \\
 \langle \mathbb{L}^{-1} \mathbf{f}, \mathbf{e}_2 \rangle 
\end{array}
\right]  
\\
&& - \dfrac{1}{2\lambda_2(k,\varepsilon)}\bigg[  \mathbb{L}^{-1}  \mathbf{e}_1 \quad  \mathbb{L}^{-1}  \mathbf{e}_2 \bigg]  
\left[
\begin{array}{cc}
\beta  &   \tilde \beta \\
\tilde \beta &\beta 
\end{array}
\right] 
\left[
\begin{array}{cc}
1 & -1 \\
-1 & 1 
\end{array}
\right] 
\left[
\begin{array}{llll}
 \langle \mathbb{L}^{-1} \mathbf{f}, \mathbf{e}_1 \rangle   \\
 \langle \mathbb{L}^{-1} \mathbf{f}, \mathbf{e}_2 \rangle 
\end{array}
\right].
\end{eqnarray*}
A further calculation yields
\begin{eqnarray*}
\boldsymbol{\varphi} &=& \mathbb{L}^{-1}  \mathbf{f}
+ \dfrac{1- \lambda_1/\langle\mathbb{L}^{-1}\mathbf{e}_1, \mathbf{e}_1+ \mathbf{e}_2 \rangle}{2\lambda_1} 
\langle \mathbb{L}^{-1} \mathbf{f}, \mathbf{e}_1 + \mathbf{e}_2  \rangle \cdot 
(\mathbb{L}^{-1}  \mathbf{e}_1 + \mathbb{L}^{-1}  \mathbf{e}_2) \\
&& + \dfrac{1- \lambda_2(k,\varepsilon)/\langle\mathbb{L}^{-1}\mathbf{e}_1, \mathbf{e}_1+ \mathbf{e}_2 \rangle}{2\lambda_2(k,\varepsilon)} 
\langle \mathbb{L}^{-1} \mathbf{f}, \mathbf{e}_1 - \mathbf{e}_2 \rangle \cdot
(\mathbb{L}^{-1}  \mathbf{e}_1 - \mathbb{L}^{-1}  \mathbf{e}_2).
\end{eqnarray*}

On the other hand, it is easy to check that 
$$
\mathbf{f}= \frac{1}{\varepsilon} 2 \cdot \mathbf{e}_1 + \kappa\cdot O(1)\cdot \mathbf{e}_1, \quad \mbox{in}\,\, V_2\times V_2. 
$$
Thus we have from Lemma \ref{lem-L_inv} that
$$
\mathbb{L}^{-1}  \mathbf{f} = \frac{1}{\varepsilon}  \left[ 2 + \kappa \cdot O(\varepsilon)\right ] \left[K^{-1}1 \cdot \mathbf{e}_1  + O(\varepsilon^{1-\tau}) \right] . 
$$
Combined with Lemma \ref{lem-identity1}, it follows that
\begin{eqnarray*}
\varepsilon \boldsymbol{\varphi} &=&  \left[ 2 + \kappa \cdot O(\varepsilon)\right ] K^{-1}1 \cdot \mathbf{e}_1  + O(\varepsilon^{1-\tau}) \\
&&+ \dfrac{1- \lambda_1/(\alpha+ O(\varepsilon^{1-\tau}))}{2\lambda_1} \cdot \left[ 2 \alpha +  O(\varepsilon^{1-\tau}) \right]  \cdot \left[ K^{-1}1 \cdot (\mathbf{e}_1 + \mathbf{e}_2) + O(\varepsilon^{1-\tau})\right] \\
&& +
 \dfrac{1- \lambda_2/(\alpha+ O(\varepsilon^{1-\tau}) )}{2\lambda_2} \cdot \left[ 2 \alpha +  O(\varepsilon^{1-\tau}) \right]  \cdot \left[ K^{-1}1 \cdot (\mathbf{e}_1 - \mathbf{e}_2) + O(\varepsilon^{1-\tau})\right]\\
&=& 
 \kappa \cdot O(\varepsilon) \cdot K^{-1}1 \cdot \mathbf{e}_1 +
\dfrac{\alpha}{\lambda_1} 
\left[K^{-1}1 \cdot (\mathbf{e}_1 + \mathbf{e}_2)  + O(\varepsilon^{1-\tau}) \right] \\
&& + \dfrac{\alpha}{\lambda_2} 
\left[K^{-1}1 \cdot (\mathbf{e}_1 - \mathbf{e}_2) + O(\varepsilon^{1-\tau}) \right] + O(\varepsilon^{1-\tau}).
\end{eqnarray*}
Similarly, we can deduce \eqref{eq-varphi-1}. This completes the proof of the lemma.  \qed \\

\subsection{ Asymptotic expansion of the solution to the scattering problem}
Define the far-field zones $\Omega_1^+:=\{ x \;|\;  x_2 > 2 \}$ and $\Omega^-_{1}:=\{ x \;|\;  x_2 <-1 \}$ above and below the slits
respectively.
\begin{lem}\label{lem-u_far_field}
The scattered field
$$ u^s_\varepsilon(x) = -\varepsilon \alpha  g^e(x,(0, 1))\cdot \left(\dfrac{1}{p} \pm \dfrac{1}{q}\right)+ O(\varepsilon^{2-\tau}) \cdot \left(\dfrac{1}{p} \pm \dfrac{1}{q}\right)$$
in $\Omega_1^\pm$.
\end{lem}

\noindent\textbf{Proof} In $\Omega_1^+$, the scattered field
$$
u^s_\varepsilon(x) = \int_{\Gamma^+_\varepsilon} g^e(x,y) \dfrac{\partial u_\varepsilon(y)}{\partial \nu} ds_y.
$$
Recall that 
$$
 \dfrac{\partial u_\varepsilon}{\partial \nu}(x_1, 1)= -\varphi_1\left(\frac{x_1}{\varepsilon}\right),
$$
we have
$$
u^s_\varepsilon(x)= -\int_{\Gamma^+_\varepsilon} g^e(x,(y_1, 1)) \varphi_1\left(\frac{y_1}{\varepsilon}\right) dy_1
= -\varepsilon \int_{0}^1 g^e(x,(\varepsilon Y, 1))  \varphi_1(Y)d Y.
$$
By  noting that
$$
g^e(x,(\varepsilon Y, 1)) = g^e(x,(0, 1))\left(1 + O(\varepsilon)\right) \quad \mbox{for} \; x\in \Omega^+ _{1},
$$
and using the asymptotic expansion for $\langle  \boldsymbol{\varphi}, \mathbf{e}_1 \rangle$ in Lemma \ref{lem-phi},
we arrive at the desired formula. The scattered field in $\Omega_1^-$ can be obtained analogously.  \qed \\

Next we consider the wave field in the slits. Since $u_\varepsilon$ is quasi-periodic, we restrict the discussion to the reference slit $S_{\varepsilon}^{(0)}$ only.  Observe that $u_\varepsilon$ satisfies
\begin{equation*} 
\left\{
\begin{array}{llll}
\vspace*{0.1cm}
\Delta u_{\varepsilon} + k^2 u_{\varepsilon} = 0 \quad \mbox{in} \; S_{\varepsilon}^{(0)},  \\
\vspace*{0.1cm}
\dfrac{\partial u_{\varepsilon}}{\partial x_1} = 0  \quad \mbox{on} \; x_1=0, \,\, x_1=\varepsilon.
\end{array}
\right.
\end{equation*}
If $k\varepsilon \ll 1$, in light of the boundary condition on the slit walls, we may expand $u_\varepsilon$ as the sum of wave-guide modes as follows:
\begin{equation}\label{u_in_slit}
u_\varepsilon(x)=a_0 \cos kx_2 + b_0 \cos k(1-x_2) + \sum_{m\geq 1} \left(a_me^{-k_2^{(m)}x_2} +  b_m e^{-k_2^{(m)}(1-x_2)} \right) \cos \frac{m\pi x_1}{\varepsilon} 
\end{equation}
where $k_2^{(m)}=\sqrt{(m\pi/\varepsilon)^2-k^2}$.

 \begin{lem}\label{lem-u_slit}
 The wave field in the slit region $S_{\varepsilon}^{(0), int} :=\{ x\in S_{\varepsilon}^{(0)} \;|\; x_2 \gg \varepsilon , 1- x_2 \gg \varepsilon \} $ is given by
$$
u_\varepsilon(x) = \bigg[ \alpha+ O(\varepsilon^{1-\tau}) \bigg] \left[ \dfrac{\cos (kx_2)}{k\sin k}  \left(\dfrac{1}{p} + \dfrac{1}{q}\right)
+ \dfrac{\cos (k(1-x_2)) }{k\sin k }   \left(\dfrac{1}{p} - \dfrac{1}{q}\right) \right] + O\left(e^{-1/\varepsilon} \right) .
$$
 \end{lem}
 \noindent\textbf{Proof}
Taking the derivative of \eqref{u_in_slit} and evaluating on the slit apertures, one has
\begin{eqnarray}\label{eqn_du_expansion}
\dfrac{\partial u_\varepsilon}{\partial x_2}(x_1,1)&=&
- a_0 k\sin k + \sum_{m\geq 1} \left( -a_m e^{-k_2^{(m)}} + b_m 
\right) k_2^{(m)} \cos \frac{m\pi x_1}{\varepsilon}, \label{eqn_du_expansion1}  \\
\dfrac{\partial u_\varepsilon}{\partial x_2}(x_1,0)&=&
b_0 k \sin k + \sum_{m\geq 1} \left( -a_m  + b_m e^{-k_2^{(m)}}
\right) k_2^{(m)} \cos \frac{m\pi x_1}{\varepsilon} \label{eqn_du_expansion2}.
\end{eqnarray}
Therefore,
\begin{eqnarray*}
- a_0 k\sin k  &=& \dfrac{1}{\varepsilon} \int_{\Gamma^+_\varepsilon} \dfrac{\partial u_\varepsilon}{\partial x_2}(x_1,1) dx_1 
= -\int_{0}^1 \varphi_1(X)d X = -\bigg[\alpha+ O(\varepsilon^{1-\tau})  \bigg] \left(\dfrac{1}{p} + \dfrac{1}{q}\right), \\
b_0 k \sin k &=& \dfrac{1}{\varepsilon} \int_{\Gamma^-_\varepsilon} \dfrac{\partial u_\varepsilon}{\partial x_2}(x_1,0) dx_1 
= \int_{0}^1 \varphi_2(X)d X = \bigg[\alpha+ O(\varepsilon^{1-\tau}) \bigg] \left(\dfrac{1}{p} - \dfrac{1}{q}\right).
 \end{eqnarray*}
Consequently
\begin{equation}\label{a0b0}
a_0  = \dfrac{1}{k \sin k}\bigg[\alpha+ O(\varepsilon^{1-\tau})  \bigg] \left(\dfrac{1}{p} + \dfrac{1}{q}\right),  \quad
b_0 = \dfrac{1}{k \sin k} \bigg[ \alpha+ O(\varepsilon^{1-\tau}) \bigg]  \left(\dfrac{1}{p} - \dfrac{1}{q}\right). 
\end{equation}

For $m\ge1$, the coefficients $a_m$ and $b_m$ can be obtained similarly by taking the inner product of \eqref{eqn_du_expansion1}
and \eqref{eqn_du_expansion2} with $\cos \dfrac{m\pi x_1}{\varepsilon}$. Then a direct estimate leads to
\begin{equation}\label{ambm}
\abs{a_m} \le C/\sqrt{m} , \quad \abs{b_m} \le C/\sqrt{m}, \quad \mbox{for}\; \,m\ge1,
\end{equation}
where $C$ is some positive constant independent of $\varepsilon$, $k$ and $m$.
The proof is complete by substituting \eqref{a0b0} and \eqref{ambm} into \eqref{u_in_slit}. \qed \\

To obtain wave field on the apertures of the slits, we consider the two apertures $\Gamma^+_\varepsilon$ and $\Gamma^-_{\varepsilon}$ of the reference slit $S_\varepsilon^{(0)}$. Define
\begin{equation} \label{eq-h}  
h(X) = \dfrac{1}{\pi}\int_{0}^1 \ln | X-Y| (K^{-1}1)(Y)dY.
\end{equation}
Let us rewrite $\beta_e(k,\kappa, d, \varepsilon)=\dfrac{1}{\pi} \ln \varepsilon +\bar \beta_e(k,\kappa,d)$, where
\begin{equation}\label{beta_bar}
\bar \beta_e(k,\kappa,d) := \dfrac{1}{\pi} \left(\ln 2 + \ln\dfrac{\pi}{d}\right) + \left(\dfrac{1}{2\pi}\sum_{n\neq 0} \dfrac{1}{|n|} - \dfrac{i}{d} \sum_{n=-\infty}^{\infty}  \dfrac{1}{\zeta_n(k)}\right).
\end{equation}

\begin{lem}\label{lem-u_aperture}
The following expansions hold for the total field
\begin{eqnarray}\label{u_up_aperture}
u_\varepsilon(x_1, 1)  &=& -\dfrac{1}{\pi}  \left( 
\dfrac{\alpha}{p} +\dfrac{\alpha}{q}  \right) \cdot  \varepsilon \ln\varepsilon - 
\left( \dfrac{\alpha}{p} + \dfrac{\alpha}{q}  \right) \left(\bar\beta_e +h(x_1/\varepsilon) \right) \cdot  \varepsilon  + 2  \nonumber \\
&& - \left(\dfrac{\alpha}{p} + \dfrac{\alpha}{q}\right) \cdot O(\varepsilon^{2-\tau}\ln\varepsilon)-\kappa \cdot O(\varepsilon)+O(\varepsilon^{2-\tau}) \nonumber
\end{eqnarray}
and
\begin{eqnarray}\label{u_low_aperture}
u_\varepsilon(x_1, 0)  &=& -\dfrac{1}{\pi}  \left( 
\dfrac{\alpha}{p} -\dfrac{\alpha}{q}  \right) \cdot  \varepsilon \ln\varepsilon - 
\left( \dfrac{\alpha}{p} - \dfrac{\alpha}{q}  \right) \left(\bar\beta_e +h(x_1/\varepsilon) \right) \cdot  \varepsilon  \nonumber \\
&& - \left(\dfrac{\alpha}{p} - \dfrac{\alpha}{q}\right) \cdot O(\varepsilon^{2-\tau}\ln\varepsilon)-\kappa \cdot O(\varepsilon)+O(\varepsilon^{2-\tau}) \nonumber
\end{eqnarray}
on the slit apertures $\Gamma^+$ and $\Gamma^-$ respectively. 
\end{lem}

 \noindent\textbf{Proof} Recall that on $\Gamma^+_\varepsilon$,
$$ u_\varepsilon(x) = \int_{\Gamma^+_\varepsilon} g_\varepsilon^e(x,y) \dfrac{\partial u_\varepsilon (y)}{\partial \nu} ds_y + u^i+ u^r.$$
Let $x_1= \varepsilon X $, $y_1= \varepsilon Y$.
We have 
\begin{equation*}
u_\varepsilon(\varepsilon X, 1) =  -\int_{0}^1 G_\varepsilon^e(X,Y)\varepsilon \varphi_1(Y)dY +f(X). 
\end{equation*} 
Using Lemma \ref{lem-phi} and the asymptotic expansion of $G_\varepsilon^e(X,Y)$ in Lemma \ref{lem-green}, we obtain
\begin{eqnarray*}
u_\varepsilon(\varepsilon X, 1) &=& -\varepsilon \beta_e\bigg(\alpha+ O(\varepsilon^{1-\tau}) \bigg) \left(\dfrac{1}{p} + \dfrac{1}{q}\right) 
 - \dfrac{\varepsilon}{\pi} \left( \kappa \cdot O(1) + \dfrac{\alpha}{p} + \dfrac{\alpha}{q} \right)\int_{0}^1 \ln | X-Y| (K^{-1}1)(Y)dY \\
&& - \left(\dfrac{\alpha}{p} + \dfrac{\alpha}{q} \right) O(\varepsilon^{2-\tau})  +  O(\varepsilon^{2-\tau}) + f(X).
 \end{eqnarray*}
The desired expansion follows by noting \eqref{eq-h} and \eqref{beta_bar}. The wave field on the lower aperture can be obtained similarly. \qed

\subsection{An overview of field enhancement and diffraction anomalies}
From Lemma \ref{lem-u_far_field} to \ref{lem-u_aperture}, we observe that $p(k;\kappa,d,\varepsilon)$ and $q(k;\kappa,d,\varepsilon)$
are two key scalar functions that will contribute to anomalous behaviors of the diffracted wave field. For instance,
the wave field will exhibit extraordinary enhancement for vanishing $p$ or $q$.
In addition, $\zeta_n(k)$, and consequently the function $\gamma(k,\kappa,d)$ (see (\ref{eq-gamma})-(\ref{eq-formula_q})) has a branch cut at certain frequencies, and this may give rise to anomalous diffracted field too. 
The following is a summary of several cases that we will explore in the rest of the paper.

\begin{enumerate}

\item [(i)]
$p(k; \kappa, d, \varepsilon)=0$ or $q(k; \kappa, d, \varepsilon)=0$ attain complex roots $k$ with negative imaginary part and real part $\mbox{Re}\;k>|\kappa|$.
Such $k$ are called resonances and the corresponding modes are called quasi-modes or leaky modes.
If the incident frequency coincides with the resonant frequency, then field enhancement will occur.

\item [(ii)]
$p(k; \kappa, d, \varepsilon)=0$ or $q(k; \kappa, d, \varepsilon)=0$ attain real roots $k$ with $k < |\kappa|$. 
Such $k$ are called real eigenvalues of the scattering operator, and the corresponding eigenmodes are called
Rayleigh-Bloch surface bound states that are confined near the periodic structure. The surface bound-state modes can be excited by nearby sources through near field interaction, but not by a plane incident wave as we consider here. 

\item [(iii)]
$p(k; \kappa, d, \varepsilon)=0$ or $q(k; \kappa, d, \varepsilon)=0$ attain real roots $k$ with $k > |\kappa|$.
In such scenario, the periodic slab structures possesses certain finite bound state embedded in the continuum states 
(or a point spectrum embedded in the continuous spectrum).

\item [(iv)]
The function $\gamma=\gamma(k, \kappa, d)$ has a branch cut at the triplet $(k, \kappa, d)$ such that $k=|\kappa + 2 \pi n/d|$ and $\zeta_n(k)=0$.
This corresponds to the Rayleigh anomaly, where the propagating mode $e^{i \kappa_n x_1 \pm i\zeta_n  x_2 }$ become an evanescent mode or vice versa.

\end{enumerate}    
We investigate (i) (ii) in Section \ref{sec-eig_res}, and explore field enhancement in Section \ref{sec-res_Rayleigh} when the resonance frequency is
close to the Rayleigh anomaly. The embedded eigenvalues (iii) is discussed briefly in Section \ref{sec-emb_eig}.

\setcounter{equation}{0}

\section{Resonances and eigenvalues away from Rayleigh cut-off frequencies}\label{sec-eig_res}
\subsection{The homogenous scattering problem}
In order to obtain eigenvalues or resonances of the scattering problem,
we consider the corresponding homogeneous problem when the incident wave $u^i=0$.
The solution $k$ is either real-valued or complex-valued with negative imaginary part. The former is called an eigenvalue and the latter is called a resonance. Here we focus on eigenvalues or resonances sufficiently away from the Rayleigh cut-off frequencies
by assuming that $\tau=0$ in \eqref{eq-delta} such that $\delta:=O(\varepsilon^{2\tau})=O(1)$.
In addition, it is natural to assume that $\delta<\pi/d$ so that $ \mathbf{R}\backslash \mathbf{B}_{\kappa,\delta}\neq\emptyset$. 
On the other hand, we only consider eigenvalues/resonances not in the high frequency regime. Therefore, we restrict to the domain of interest to
$$ D_{\kappa,\delta, M} := \{ z: \;\; |z| \le M\} \backslash B_{\kappa,\delta},$$
where $M>0$ is a fixed constant.

As to be shown later on, for each $\kappa\in(-\pi/d,\pi/d]$, the eigenvalues and resonances in $D_{\kappa,\delta, M}$ lie in the vicinity of 
$$k_{m,0}:=m\pi \quad \mbox{for} \; m=1, 2, \cdots,  \; \mbox{and} \; m\pi < M.$$ 
Therefore, we extend the asymptotic expansions of the boundary integral operators to a neighborhood of $k_{m,0}$ on the complex plane. 
More precisely, let $\delta_0=\min\{\delta/2, \pi\}$. If $k_{m,0}\in D_{\kappa,\delta, M}$, we choose the disc with radius $\delta_0$ centered at $k_{m,0}$
on the complex $k$-plane, which is denoted as $B_{\delta_0}(k_{m,0})$. We analytically extend the functions $\beta_e(k)$,  $\beta_i(k)$,  and $\tilde \beta(k)$,  which are defined in \eqref{beta_e} - \eqref{t_beta} for real $k$, to  the neighborhood of $k_{m,0}$, $B_{\delta_0}(k_{m,0})$.  
One can show that the asymptotic expansions for the kernels $G_\varepsilon^e$,  $G_\varepsilon^i$,  and $\tilde G_\varepsilon^i$ given in Lemma \ref{lem-green} holds in  $B_{\delta_0}(k_{m,0})$.

By virtue of \eqref{eq-scattering4}, the homogeneous problem is equivalent to the operator equation
$$(\mathbb{P} + \mathbb{L})\boldsymbol{\varphi}=0 \quad \mbox{in} \; B_{\delta_0}(k_{m,0}). $$
In light of \eqref{eqn-linear_sys}, this reduces to
\begin{equation*}
( \mathbb{M}+\mathbb{I})\left[
 \begin{array}{llll}
\langle  \boldsymbol{\varphi}, \mathbf{e}_1 \rangle  \\
\langle  \boldsymbol{\varphi}, \mathbf{e}_2 \rangle
\end{array}
\right] =0,
\end{equation*}
where the matrix $\mathbb{M}$ is defined by \eqref{eqn-matrix-m}.
Note that the eigenvalues of $\mathbb{M}+\mathbb{I}$ are given by \eqref{eq-lambda1} and \eqref{eq-lambda2},
thus the characteristic values of the operator-valued function $\mathbb{P} + \mathbb{L}$ are the roots of the two analytic functions $\lambda_1(k; \kappa, d, \varepsilon) $ and $\lambda_1(k; \kappa, d, \varepsilon) $, or equivalently,  $p(k; \kappa, d, \varepsilon) $ and $q(k; \kappa, d, \varepsilon) $
as defined in \eqref{eq-def_pq}.

\begin{lem} \label{lem-equi}
The resonances of the scattering problem \eqref{eq-Helmholtz} - \eqref{eq-rad_cond} in $B_{\delta_0}(k_{m,0})$ are the roots of one of the analytic functions $p(k; \kappa, d, \varepsilon) =0$ and $q(k; \kappa, d, \varepsilon) =0$ with $\Im \,k<0$, and the eigenvalues are those with $\Im \,k =0$.
\end{lem}

\begin{lem}\label{lem-asym_res_eig}
For each $\kappa\in(-\pi/d,\pi/d]$, the roots of $p(k; \kappa, d, \varepsilon) =0$ and $q(k; \kappa, d, \varepsilon) =0$
in the domain $D_{\kappa, \delta, M}$ attain the following asymptotic expansion:
\begin{equation}\label{eqn-asym_res_eig}
k_m= k_m(\kappa, d, \varepsilon)=  m\pi + 2m \pi  \left [ \dfrac{1}{\pi}\varepsilon \ln\varepsilon + \left(  \dfrac{1}{\alpha} +  \gamma(m\pi, \kappa, d)  \right)\varepsilon \right] + O(\varepsilon^2\ln^2\varepsilon)
\end{equation}
for even and odd integers $m$ respectively. 
\end{lem}
\noindent\textbf{Proof}. 
Let $c(k) =\dfrac{\cot k }{k} + \dfrac{1}{k\sin k}$. From \eqref{eq-formula_p}, it follows that
\begin{equation}\label{res_lambda1}
p(k; \kappa, d, \varepsilon)=\varepsilon +\left[ c(k) + \varepsilon \gamma(k,\kappa,d) + \dfrac{1}{\pi} \varepsilon \ln \varepsilon   \right]   \left(\alpha + r(k, \varepsilon) \right) =0,
\end{equation}
where $\gamma(k,\kappa,d)$ is defined in \eqref{eq-gamma}, $r(k, \varepsilon)$ is analytic in $B_{\delta_0}(k_{m,0})$ and  $r(k, \varepsilon) \sim O(\varepsilon)$.  It is clear that the analytic function $c(k)$ attains a simple root $k_{m,0}=m\pi$ in $B_{\delta_0}(k_{m,0})$ for odd integers $m$.
From Rouche's theorem, we deduce that there is a simple root  of $p(k,\varepsilon)$, which is denoted as $k_{m}$,  close to $k_{m,0}$
if $\varepsilon$ is sufficiently small. 

To obtain the leading-order asymptotic terms of $k_{m}$, first let us consider the root for
\begin{equation}\label{p1}
p_1(k; \kappa, d, \varepsilon):=\varepsilon +\left[ c(k) + \varepsilon \gamma(k,\kappa,d) + \dfrac{1}{\pi} \varepsilon \ln \varepsilon   \right] \alpha=0.
\end{equation}
The Taylor expansion for $p_1(k, \varepsilon)$ at $k=k_{m,0}$ yields
\begin{eqnarray}\label{p1_expansion}
p_1(k; \kappa, d, \varepsilon) &=& \varepsilon +\bigg[ c'(k_{m,0})(k-k_{m,0}) +O(k-k_{m,0})^2 + \varepsilon \gamma(k_{m,0}) \\ 
&& +\varepsilon \cdot O(k-k_{m,0})  + \dfrac{1}{\pi} \varepsilon \ln \varepsilon   \bigg]  \alpha. \nonumber
\end{eqnarray}
A direct calculation gives $c'(k_{m,0})= -\dfrac{1}{2m\pi}$. We can deduce that $p_1$ has a simple root $k_{m,1}$ close to $k_{m,0}$, and is given by 
$$
k_{m, 1} =k_{m, 0}+ 2m\pi \left [ \dfrac{1}{\pi}\varepsilon \ln\varepsilon + \left(  \dfrac{1}{\alpha} + \gamma_0(k_{m,0},\kappa,d) \right)\varepsilon \right] + O(\varepsilon^2\ln^2\varepsilon). 
$$

Next we show that $k_{m}- k_{m, 1} = O(\varepsilon^2\ln^2\varepsilon )$ and the desired asymptotic expansion \eqref{eqn-asym_resonance} follows.
Note that 
$$
p(k; \kappa, d, \varepsilon) - p_1(k; \kappa, d, \varepsilon) = O(c(k)+\varepsilon\ln\varepsilon) \; r(k, \varepsilon)
$$
and
$$  
p_1(k; \kappa, d, \varepsilon) = c(k)\alpha + O(\varepsilon\ln\varepsilon).
$$
Hence, one can find a constant $C_m>0$ such that 
$$
| p(k; \kappa, d, \varepsilon) - p_1(k; \kappa, d, \varepsilon)|  < | p_1(k; \kappa, d, \varepsilon)|
$$ 
for all $k$ such that $|k-k_{m, 1}| = C_m \varepsilon^2\ln^2\varepsilon$. 
 By the Rouche's theorem,  $p$ has a simple root in the disc  $\{k \;|\; \;|k-k_{m, 1}| = C_m \varepsilon^2\ln^2\varepsilon\}$, which proves our claim. 

Similarly, we obtain the root of $q(k; \kappa, d, \varepsilon) =0$ in $B_{\delta_0}(k_{m,0})$ for even integers $m$.
The arguments are the same as above and we omit here.  \qed

\subsection{Asymptotic expansions of resonances and eigenvalues}\label{sec-asy_res_eig}
The formula (\ref{eqn-asym_res_eig}) gives the asymptotic expansion of resonances and eigenvalues for the scattering problem \eqref{eq-Helmholtz} - \eqref{eq-rad_cond}. We now distinguish between resonances and eigenvalues.
First, observe that  
$$\Im\,\gamma(m\pi, \kappa, d)\neq 0  \quad \mbox{if} \; m\pi>|\kappa|.$$
Therefore, $k_m$ attains a non-zero imaginary part, and $k_m$ is a complex-valued resonance.
We immediately have the following proposition.
\begin{prop}\label{prop-resonance}
For each $\kappa\in(-\pi/d,\pi/d]$, if $m\pi>|\kappa|$ such that $\Im \, \gamma(m\pi, \kappa, d) \neq 0$,
then there exists resonance in the domain $D_{\kappa,\delta, M}$ for the scattering problem \eqref{eq-Helmholtz} - \eqref{eq-rad_cond} with the following asymptotic expansion
\begin{equation}\label{eqn-asym_resonance}
k_m =  m\pi + 2m \pi  \left [ \dfrac{1}{\pi}\varepsilon \ln\varepsilon + \left(  \dfrac{1}{\alpha} +  \gamma(m\pi, \kappa, d)  \right)\varepsilon \right] + O(\varepsilon^2\ln^2\varepsilon).
\end{equation}
 Here $\alpha=\langle K^{-1} 1, 1 \rangle$.
\end{prop}

\begin{rmk}
In Proposition \ref{prop-resonance}, the imaginary part of the resonance $k_{m}$ is given by $2m\pi \cdot \mbox{Im} \gamma(m\pi, \kappa, d) \cdot \varepsilon + O(\varepsilon^2\ln^2\varepsilon)$, which is negative and is of order $O(\varepsilon)$.
As the size of period $d\to+\infty$, we see that
$$ \mbox{Im} \gamma(m\pi,\kappa, d) = -\dfrac{1}{d} \sum_{|\kappa_n|<m\pi} \dfrac{1}{\zeta_n(m\pi) } \to -\dfrac{1}{2\pi}\int_{-m\pi}^{m\pi} \dfrac{1}{\sqrt{(m\pi)^2-t^2}}dt = -\dfrac{1}{2}. $$
This constant is consistent with the one obtained for the case of a single slit perforated in a slab of infinite length (cf \cite{lin_zhang16}). \\
\end{rmk}

Now if $m\pi<|\kappa|$ such that 
$$\Im \gamma(m\pi, \kappa, d) = 0,$$ 
the imaginary part of $O(\varepsilon)$-term of $k_m$ in \eqref{eqn-asym_res_eig} is zero. 
However, one can not tell directly from the asymptotic expansion (\ref{eqn-asym_res_eig})
whether the higher order terms of $k_m$ are real or complex-valued.
Instead, we resort to the variational formulation to verify that those $k_m$ are indeed real eigenvalues.

Let us first recall some basic facts from \cite{bonnet_starling94}.
Denote $\Omega_{\varepsilon}^{(0)}= \Omega^{(0)} \cap \Omega_{\varepsilon}$ and define the function space
$$
H^1_{\kappa,d }(\Omega_{\varepsilon}^{(0)})= \Big\{u: u \in H^1(\Omega_{\varepsilon}^{(0)}), u(d, x_2)= e^{i\kappa d} u(0, x_2), \frac{\partial u}{\partial x_1}(d, x_2) = e^{i\kappa d}\frac{\partial u}{\partial x_1}(0, x_2) \Big\}.
$$
We define a sesquilinear form
$$
a(u, v) := \int_{\Omega_{\varepsilon}^{(0)}} \nabla u \cdot \nabla \bar v dx
$$
on $H^1_{\kappa,d}(\Omega_{\varepsilon}^{(0)}) \times H^1_{\kappa,d}(\Omega_{\varepsilon}^{(0)})$ and denote 
$$
(u, v)= \int_{\Omega_{\varepsilon}^{(0)}} u \bar v dx.
$$ 
Let $A(\kappa)=A(\kappa, d, \varepsilon)$ be the operator associated with the sesquilinear form $a$ such that
$$
(A(\kappa)u, v) = a(u, v)
$$
for all $u, v \in H^1_{\kappa,d}(\Omega_{\varepsilon}^{(0)})$. 
The following statement holds (cf. \cite{bonnet_starling94}).

\begin{lem}\label{lem-opA}
\begin{enumerate}
\item
$(k,\kappa, u)$ is a solution to the homogeneous problem with $u^i=0$ if and only if $k^2$ is an eigenvalue of $A(\kappa)$ and $u$ is the associated eigenfuction. 
\item
$A(\kappa)$ is a positive self-adjoint operator. 
\end{enumerate}
\end{lem} 

For each positive integer $m$, let
$$
\Lambda_m(\kappa) = \Lambda_m(\kappa, d, \varepsilon)= \inf_{ V_m \in \mathcal{V}_m}  \sup_{ u\in  V_m, u \neq 0} \frac{a(u,u)}{(u, u)},
$$
where $\mathcal{V}_m$ denotes the set of all m-dimensional subspaces of $H^1_{\kappa,d}(\Omega_{\varepsilon})$. 
It is clear that for each fixed $\kappa$, $\Lambda_m(\kappa)$ is an increasing sequence. 
We denote by $\mathcal{N}(\kappa)= \mathcal{N}(\kappa, d, \varepsilon)$ the number of eigenvalues $\Lambda$ of $A(\kappa)$, counting their multiplicity, 
which are strictly less than $|\kappa|^2$.

\begin{lem}\label{lem-Lambda} The following statements hold (cf. \cite{bonnet_starling94}):
\begin{enumerate}

\item
$\Lambda_m(\kappa) = \Lambda_m(-\kappa)$ for $\kappa \in [0, \pi/d]$.

\item
For each fixed $m$, $\Lambda_m(\kappa)$ is a continuous function of $\kappa$. 

\item For $\kappa\in(-\pi/d, \pi/d]$, 
if $k_{m,0}<|\kappa|$, then  $\Lambda_m(\kappa)<\kappa^2$ and $\mathcal{N}(\kappa)\ge m$. 
In addition, $\Lambda_1$, $\Lambda_2$, $\cdots$, $\Lambda_m$
are the first $m$ eigenvalues of $A(\kappa)$. 

\end{enumerate}

\end{lem}

\begin{figure}[!htbp]
\begin{center}
\includegraphics[height=4.2cm,width=12cm]{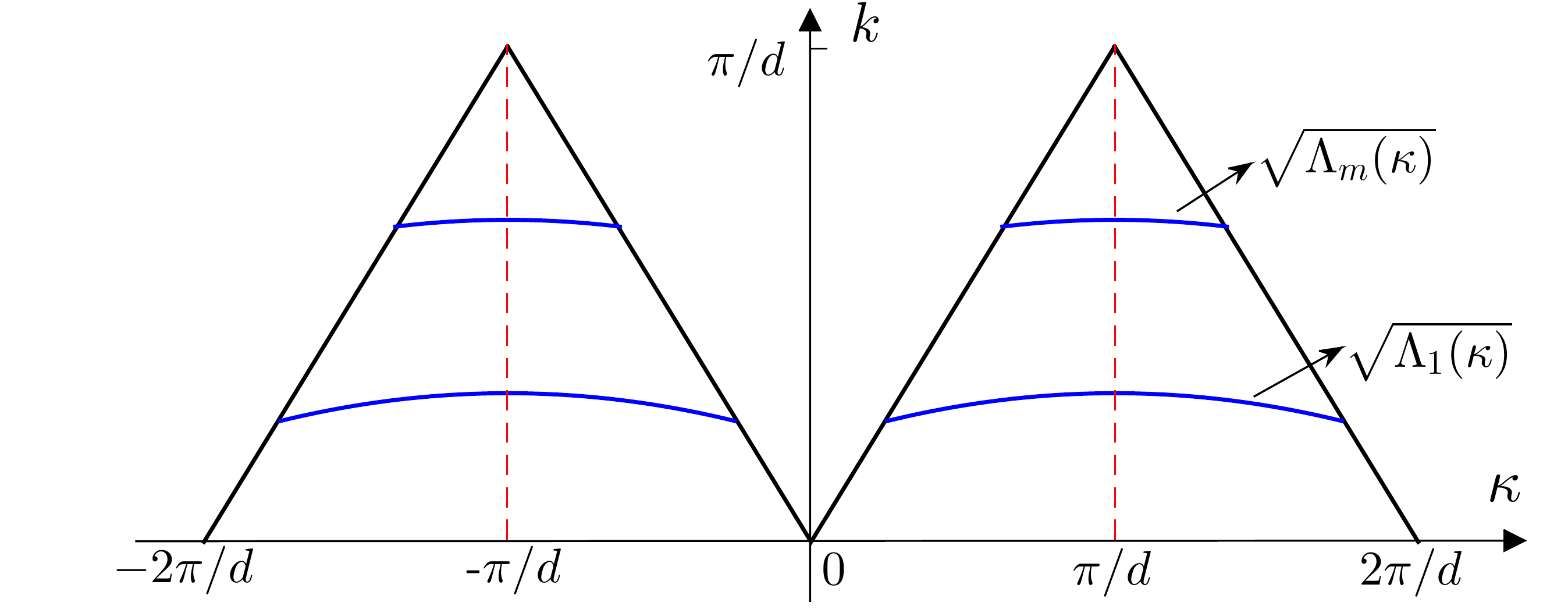}
\vspace{-10pt}
\caption{Eigenvalues of the scattering operator $A(\kappa)$.}\label{fig-eigenvalue}
\vspace{-10pt}
\end{center}
\end{figure}

Now let us first consider the special case when $\kappa=\pi/d$ and $m\pi<|\kappa|$. 
Then for real $k$ in the neighborhood of $k_{m,0}$, the Green's function $G_\varepsilon^e(X,Y; \pi/d )$ is a real-valued function,
by noting that the $n$ and $-(n+1)$ terms in the following series expansion 
$$ G_\varepsilon^e(X,Y; \pi/d) = -\dfrac{i}{d}\sum_{n=-\infty}^{\infty} \dfrac{1}{\zeta_n(k) } e^{i \kappa_n\varepsilon(X-Y)} =  -\dfrac{1}{d}\sum_{n=-\infty}^{\infty} \dfrac{1}{ \sqrt{((2n+1)\pi/d)^2 - k^2}} e^{i (2n+1)\pi/d\varepsilon(X-Y)} $$
are conjugate pairs.
Therefore, all the terms in \eqref{eq-lambda1} and \eqref{eq-lambda2} are real. Consequently, $p(k; \pi/d, d, \varepsilon)=0$ or $q(k; \pi/d, d, \varepsilon)=0$ attains a real root $k_m$ near $m\pi$, which is an eigenvalue of the scattering problem. The 
asymptotic expansion of the eigenvalue is given by  (\ref{eqn-asym_res_eig}).
In particular, by Lemma \ref{lem-opA} and \ref{lem-Lambda}, we have 
$$
\Lambda_{m}(\pi/d, d, \varepsilon) = k_m^2(\pi/d, d, \varepsilon).  
$$
For $\kappa\neq\pi/d$, the continuity of $\Lambda_m(\kappa)$ and $k_m(\kappa)$ implies that  $\Lambda_{m}(\kappa, d, \varepsilon) = k_m^2(\kappa, d, \varepsilon)$ as long as $m\pi<|\kappa|$, and $k_m(\kappa)$ is a real eigenvalue. In summary, we can draw the following conclusion.
\begin{prop}\label{prop-eigen}
For each $\kappa\in(-\pi/d,\pi/d]$, if $m\pi<|\kappa|$ such that  $\Im \gamma(m\pi, \kappa, d) = 0$,
then there exists exactly $m$ eigenvalues in the domain $D_{\kappa,\delta, M}$ for the scattering problem \eqref{eq-Helmholtz} - \eqref{eq-rad_cond}. Moreover, the following asymptotic expansion holds for each eigenvalue $k_m$:
\begin{equation}\label{eqn-asym_resonance}
k_m =  m\pi + 2m \pi  \left [ \dfrac{1}{\pi}\varepsilon \ln\varepsilon + \left(  \dfrac{1}{\alpha} +  \gamma(m\pi, \kappa, d)  \right)\varepsilon \right] + O(\varepsilon^2\ln^2\varepsilon).
\end{equation}
Here $\alpha=\langle K^{-1} 1, 1 \rangle$.
\end{prop}

Note that $k_m<|\kappa|$ and $\zeta_n(k_m)$ is a pure imaginary number for all $n$, we see that the associated eigenfunctions are Rayleigh-Bloch surface bound states that decays exponentially away from the slab structure.

\begin{rmk} \label{rmk-11}
Note that if the size of period $d$ satisfies the inequality 
$$ k_1(\pi/d, d, \varepsilon) \geq \pi/d,$$ then $\mathcal{N}(\pi/d)=0$. Therefore, there is no eigenvalue for the scattering problem and only resonances exist. On the other hand, if $$ k_1(\pi/d, d, \varepsilon)< \pi/d,$$  then $\mathcal{N}(\pi/d, d, \varepsilon) \geq 1$.
That is, in addition to resonances, there exists at least one eigenvalue which is near $\pi$.
\end{rmk}

\subsection{Field enhancement at resonant frequencies}\label{sec-field_res}
To investigate the field enhancement,  recall that $p(k; \kappa, d, \varepsilon):=\varepsilon \lambda_1(k; \kappa, d, \varepsilon)$ and 
$q(k; \kappa, d, \varepsilon):=\varepsilon \lambda_2(k; \kappa, d, \varepsilon)$.

\begin{lem} \label{lem-pq}
If $m\varepsilon \ll 1$, then at the resonant frequencies $k= \Re \, k_m$,
$$ p(k; \kappa, d, \varepsilon)= i \cdot \Im \, \gamma(m\pi, \kappa, d) \cdot \alpha\varepsilon + O( \varepsilon^2 \ln^2\varepsilon) $$
and
$$ q(k; \kappa, d, \varepsilon)= i  \cdot \Im \, \gamma(m\pi, \kappa, d) \cdot \alpha \varepsilon + O( \varepsilon^2 \ln^2\varepsilon)  $$
where $m$ is odd and even respectively.
\end{lem}

\noindent\textbf{Proof}
We expand $p(k, \varepsilon)$ and $q(k, \varepsilon)$  in the disk $\{ | k- \Re \, k_{m} | \leq \varepsilon |\ln \varepsilon| \}$.
If $m$ is odd,  from the definition of $p_1$ in \eqref{p1} and its expansion \eqref{p1_expansion}, it follows that
\begin{eqnarray*}
p(k; \kappa, d, \varepsilon) &=& p_1(k; \kappa, d, \varepsilon)+  O( \varepsilon^2 \ln \varepsilon) \\
& =&  p_1'(k_{m})(k-k_{m}) + O(k-k_{m})^2 + O( \varepsilon^2 \ln \varepsilon) \\
&= & \left[ \alpha c'(k_{m,0}) + O(\varepsilon \ln \varepsilon) \right]\cdot(k-k_{m}) + O( \varepsilon^2\ln^2\varepsilon) \\
&= &  \alpha  c'(k_{m,0})\cdot (k-k_{m}) + O( \varepsilon^2 \ln^2\varepsilon) \\
&=&  -\frac{\alpha }{2m\pi} \left( k- \Re \, k_{m} - i \; \Im \,k_{m} \right)+ O( \varepsilon^2 \ln^2\varepsilon).
\end{eqnarray*}
Note that 
$$
\Im \, k_{m} = 2m\pi \cdot \Im \, \gamma(m\pi, \kappa, d) \cdot \varepsilon  +  O( \varepsilon^2 \ln^2\varepsilon).
$$
We deduce that at the odd resonant frequencies $k= \Re \, k_{m}$,
$$
p(k; \kappa, d, \varepsilon)= i\cdot \Im \, \gamma(m\pi, \kappa, d) \cdot \alpha \varepsilon + O( \varepsilon^2 \ln^2\varepsilon).
$$ 
The calculations for $q(k; \kappa, d, \varepsilon)$ at the even resonant frequencies are the same.  \qed \\

The following proposition follows directly from Lemma \ref{lem-phi} and  Lemma \ref{lem-pq} .
\begin{prop}
At resonant frequencies,  $\boldsymbol{\varphi}\sim O\left(1/\varepsilon\right)$  in $V_1 \times V_1$ and $\langle  \boldsymbol{\varphi}, \mathbf{e}_i \rangle \sim O\left(1/\varepsilon\right)$, i=1,2.
\end{prop}

We now investigate the field enhancement in far-field and near-filed zones.
From Lemma \ref{lem-u_far_field}, the scattered field
$$ u^s_\varepsilon(x) = -\varepsilon \alpha  g^e(x,(0, 1))\cdot \left(\dfrac{1}{p} \pm \dfrac{1}{q}\right)+ O(\varepsilon^{2-\tau}) \cdot \left(\dfrac{1}{p} \pm \dfrac{1}{q}\right) \quad \mbox{for} \; x\in \Omega^\pm _{1}. $$
At the resonant frequencies $k= \Re \, k_{m}$ when $m$ is odd, an application of Lemma \ref{lem-pq} yields 
$$
\dfrac{1}{p}= -\frac{i}{ \Im\,\gamma(m\pi,\kappa,d) \alpha \varepsilon}(1+  O( \varepsilon \ln^2 \varepsilon)).
$$
Hence the scattered field
\begin{equation}\label{eqn-u_far-field1}
u^s_\varepsilon(x)= \frac{i}{ \Im\,\gamma(m\pi,\kappa,d)} \cdot  g^e(x,(0, 1))+O(\varepsilon \ln^2 \varepsilon) \quad \mbox{for} \; x\in \Omega^+ _{1}.
\end{equation}
It is seen that the scattering enhancement is of order $O(\varepsilon^{-1})$ compared to the 
$O(\varepsilon)$ order for the scattered field at non-resonant frequencies. 
In addition, the scattered field behaves as the radiating field of a periodic array of monopoles located at 
$\displaystyle{\bigcup_{n=-\infty}^{\infty} (nd,1)}$.
The same occurs at resonant frequencies  $k= \Re \, k_{m}$ when $m$ is even by an application of Lemma \ref{lem-pq}. 

Following the similar calculations as above, in the far-field zone $\Omega^-_{1}$ below the slits, the transmitted field is equivalent to the radiating field of an array of monopoles located at $\displaystyle{\bigcup_{n=-\infty}^{\infty} (nd,0)}$, and is given by
 \begin{eqnarray*}
u^s_\varepsilon(x)
 = -\varepsilon \alpha  g^e(x,(0, 0))\cdot \left(\dfrac{1}{p} - \dfrac{1}{q}\right)+ O(\varepsilon^2) \cdot \left(\dfrac{1}{p} - \dfrac{1}{q}\right).
\end{eqnarray*}
It follows that 
\begin{equation}\label{eqn-u_far-field2}
u^s_\varepsilon(x)= \frac{i}{ \Im\,\gamma(m\pi,\kappa,d)} \cdot  g^e(x,(0, 0))+O(\varepsilon \ln^2 \varepsilon)
\end{equation}
and
\begin{equation}\label{eqn-u_far-field3}
u^s_\varepsilon(x)= -\frac{i}{ \Im\,\gamma(m\pi,\kappa,d)} \cdot  g^e(x,(0, 0))+O(\varepsilon \ln^2 \varepsilon) 
\end{equation}
in $\Omega^- _{1}$ at the odd and even resonant frequencies respectively. The transmission enhancement is of order $O(\varepsilon^{-1})$ at the resonant frequencies. \\

\noindent\textbf{Remark 4.1} \textit{The amplitude of the scattered and transmitted field also depends on the reduced wave vector $\kappa$ and 
the size of the period $d$. This is explicitly given by the scalar function $\Im\,\gamma(m\pi,\kappa,d)$. The same holds true
in the near field described below.} \\

The shape of resonant wave modes in the slits  and their enhancement orders are characterized in the following theorem.
 \begin{thm}\label{thm-u_slit_res}
 The wave field in the slit region $S_{\varepsilon}^{(0), int} :=\{ x\in S_{\varepsilon}^{(0)} \;|\; x_2 \gg \varepsilon , 1- x_2 \gg \varepsilon \} $ is given by
 \begin{equation}\label{eqn-u_slit1}
 u_\varepsilon(x) =  - \left(\dfrac{1}{\varepsilon} + O(\ln^2\varepsilon)  +  O(1)  \right) \cdot \dfrac{i \cdot \cos (k(x_2-1/2))}{\Im\,\gamma(m\pi,\kappa,d) \cdot k\sin (k/2)} + \dfrac{\sin (k(x_2-1/2)) }{\sin (k/2)} + O(\varepsilon \ln\varepsilon) 
\end{equation}
\begin{equation}\label{eqn-u_slit2}
u_\varepsilon(x) =  \left(\dfrac{1}{\varepsilon} + O(\ln^2\varepsilon)  + O(1)  \right)  \cdot \dfrac{i\cdot \sin (k(x_2-1/2))}{\Im\,\gamma(m\pi,\kappa,d) \cdot k\cos (k/2)} + \dfrac{\cos (k(x_2-1/2)) }{\cos (k/2)} + O(\varepsilon \ln\varepsilon)
\end{equation}
at the resonant frequencies $k= \Re \, k_{m}$ where $m$ is odd and even respectively.
 \end{thm}
\noindent\textbf{Proof} 
By Lemma  \ref{lem-u_slit}, in the region $S_{\varepsilon}^{(0),int} $,
 \begin{eqnarray*}
u_\varepsilon(x) &=&  \bigg[ \alpha+ O(\varepsilon) \bigg] \left[ \dfrac{\cos (kx_2)}{k\sin k}  \left(\dfrac{1}{p} + \dfrac{1}{q}\right)
+ \dfrac{\cos (k(1-x_2)) }{k\sin k }   \left(\dfrac{1}{p} - \dfrac{1}{q}\right) \right] + O\left(e^{-1/\varepsilon} \right) \\
&=& 2 \bigg[ \alpha+ O(\varepsilon)  \bigg] \left[ \dfrac{1}{p} \dfrac{\cos(k/2)\cos (k(x_2-1/2)) }{k\sin k}  - \dfrac{1}{q}\dfrac{\sin(k/2)\sin (k(x_2-1/2)) }{k\sin k }  \right]  \\
&& +  O\left(e^{-1/\varepsilon} \right).
\end{eqnarray*}
At resonant frequencies $k= \Re \, k_{m}$ when $m$ is odd,
\begin{equation*}
\dfrac{1}{p}= -\frac{i}{ \Im\,\gamma(m\pi,\kappa,d) \alpha \varepsilon}(1+  O( \varepsilon \ln^2 \varepsilon)) \quad \mbox{and} \quad
\dfrac{1}{q}= \dfrac{k \sin k}{ (\cos k -1) \alpha} (1 + O(\varepsilon \ln\varepsilon)).
\end{equation*}
Therefore,
\begin{eqnarray*}
u_\varepsilon(x)&=&\left(1+ O(\varepsilon) + O(\varepsilon^2\ln\varepsilon) \right) \bigg[ -\dfrac{1}{\varepsilon} \cdot \dfrac{i\cos (k(x_2-1/2)) }{\Im\,\gamma(m\pi,\kappa,d) \cdot k\sin (k/2)} (1 + O(\varepsilon \ln^2\varepsilon)) \\
&& +  \dfrac{\sin (k(x_2-1/2)) }{\sin(k/2)} (1 + O(\varepsilon \ln\varepsilon)) \bigg] +  O\left(e^{-1/\varepsilon} \right)  \\
&=& -\left(\dfrac{1}{\varepsilon} + O(\ln^2\varepsilon)  +  O(1)  \right)  \cdot \dfrac{i}{\Im\,\gamma(m\pi,\kappa,d) \cdot k\sin (k/2)} \cdot \cos (k(x_2-1/2)) \\
&&+ \dfrac{\sin (k(x_2-1/2)) }{\sin (k/2)} + O(\varepsilon \ln\varepsilon).
\end{eqnarray*}
Similarly, at resonant frequencies $k= \Re \, k_{m}$ when $m$ is even,
\begin{equation*}
\dfrac{1}{p}= \dfrac{k \sin k}{ (\cos k +1) \alpha} (1 + O(\varepsilon \ln\varepsilon)) \quad \mbox{and} \quad
\dfrac{1}{q}= -\frac{i}{ \Im\,\gamma(m\pi,\kappa,d) \alpha \varepsilon}(1+  O( \varepsilon \ln^2 \varepsilon)), 
\end{equation*}
and we obtain
\begin{eqnarray*}
u_\varepsilon(x)&=&\left(1+ O(\varepsilon) + O(\varepsilon^2\ln\varepsilon) \right)
 \bigg[\dfrac{\cos (k(x_2-1/2)) }{ k\cos (k/2)} (1 + O(\varepsilon \ln\varepsilon)) \\
&& +\dfrac{1}{\varepsilon} \cdot \dfrac{i\sin (k(x_2-1/2)) }{\Im\,\gamma(m\pi,\kappa,d) \cdot k\cos(k/2)} (1 + O(\varepsilon \ln^2\varepsilon)) \bigg] +  O\left(e^{-1/\varepsilon} \right)  \\
&=&  \left(\dfrac{1}{\varepsilon} + O(\ln^2\varepsilon)  + O(1)  \right) \cdot \dfrac{i}{\Im\,\gamma(m\pi,\kappa,d) \cdot k\cos (k/2)} \cdot \sin (k(x_2-1/2)) \\
&& + \dfrac{\cos (k(x_2-1/2)) }{\cos (k/2)} + O(\varepsilon \ln\varepsilon).
\end{eqnarray*} 
\qed

Therefore, the enhancement due to resonance is of order $O(\varepsilon^{-1})$ in the slit. 
Moreover, the dominant modes in the slit takes the simple form of $ \cos (k(x_2-1/2))$ and $ \sin (k(x_2-1/2))$ at 
odd and even resonant frequencies respectively.

An application of Lemma \ref{lem-u_aperture} and Lemma \ref{lem-pq} leads to the enhanced field on the apertures of slits as stated below.
\begin{thm}
Let function $h$ and $\bar\beta_e$ be defined by \eqref{eq-h} and \eqref{beta_bar} respectively. 
The wave fields on the apertures $\Gamma^+_\varepsilon$ and $\Gamma^-_{\varepsilon}$  are 
\begin{eqnarray}\label{eqn-slit_aperture1}
u_\varepsilon(x_1, 1) &=& \dfrac{i}{\pi \cdot \Im\,\gamma(m\pi,\kappa,d)} \cdot \left[\ln\varepsilon + \left(\bar\beta_e + h(x_1/\varepsilon) \right)\right] +2+ O(\varepsilon \ln^3 \varepsilon) \\
u_\varepsilon(x_1, 0) &=& \dfrac{i}{\pi \cdot \Im\,\gamma(m\pi,\kappa,d)} \cdot \left[\ln\varepsilon + \left(\bar\beta_e + h(x_1/\varepsilon) \right)\right]  + O(\varepsilon \ln^3 \varepsilon)
\end{eqnarray}
at the resonant frequencies $k= \Re \, k_{m}$ for odd $m$.
At the resonant frequencies $k= \Re \, k_{m}$ when $m$ is even, the wave fields on the slit apertures $\Gamma^+_\varepsilon$ and $\Gamma^-_{\varepsilon}$  are 
\begin{eqnarray}\label{eqn-slit_aperture2}
u_\varepsilon(x_1, 1) &=& \dfrac{i}{\pi \cdot \Im\,\gamma(m\pi,\kappa,d)} \cdot \left[\ln\varepsilon + \left(\bar\beta_e + h(x_1/\varepsilon) \right)\right]  +2+ O(\varepsilon \ln^3 \varepsilon) \\
u_\varepsilon(x_1, 0) &=& -\dfrac{i}{\pi \cdot \Im\,\gamma(m\pi,\kappa,d)} \cdot \left[\ln\varepsilon + \left(\bar\beta_e + h(x_1/\varepsilon) \right)\right] + O(\varepsilon \ln ^3\varepsilon). 
\end{eqnarray}
\end{thm}
It is seen that the leading-order of the resonant mode is a constant of order $O\ln (\varepsilon)$ along the slit apertures,
and the enhancement due to the resonant scattering is of order $O(\varepsilon^{-1})$.

\section{Field enhancement at resonant frequencies near Rayleigh anomaly}\label{sec-res_Rayleigh}
Rayleigh anomaly occurs at cut-off frequencies when $k= \pm (\kappa + 2\pi n/d)$ for some integer $n$ such that $\zeta_n(k)=0$. 
According to the Rayleigh-Bloch expansion \eqref{eq-rad_cond}, this corresponds to a grazing angle near which 
the propagating mode $e^{i \kappa_n x_1 \pm i\zeta_n  x_2 }$ becomes an evanescent mode and vice versa. In this section,
we investigate the field enhancement at those resonant frequencies that are close to the Rayleigh cut-off frequencies,
by assuming that $0<\tau<1$ in \eqref{eq-delta} and  $\delta=O(\varepsilon^{2\tau}) \ll 1$.

We consider a pair $(k^0, \kappa^0)$ which satisfies
$k^0= \kappa^0 + 2\pi n_0/d$ for some integer $n_0$ and $k^0$ is a cut-off frequency.
For clarity of presentation, let us assume that $\kappa^0\neq 0, \pi/d$ so that only one mode among all diffracted orders
turns into an evanescent mode near $(k^0, \kappa^0)$. 
However, the following derivations can be extended straightforwardly for  $\kappa^0 = 0, \pi/d$.
Suppose that $k$ is perturbed away from the cut-off frequency $k^0$ such that 
$kd = k^0 d+\delta = \kappa^0 d + 2\pi n_0 +\delta$, where $\delta=O(\varepsilon^{2\tau})$ and $0<\tau<1$.
It is clear that $e^{i \kappa_{n_0} x_1 \pm i\zeta_n  x_2 }$ is a propagating mode if $\delta>0$ and an evanescent mode if $\delta<0$.

A direct expansion yields
$$ \dfrac{1}{d \cdot \zeta_{n_0}(k)}=\frac{1}{\sqrt{|\delta|}} e^{-\frac{1}{2}i \arg \delta} \left( \frac{1}{\sqrt{2k_0d}}+ O(\delta)\right). $$
From the definition of $\gamma$ in \eqref{eq-gamma}, we see that
\begin{eqnarray*}
\gamma(k, \kappa^0, d) &=& \dfrac{1}{\pi} \left(3\ln 2 + \ln\dfrac{\pi}{d}\right) + \left(\dfrac{1}{2\pi}\sum_{n\neq 0} \dfrac{1}{|n|} 
- \dfrac{i}{d} \sum_{n\neq n_0}  \dfrac{1}{\zeta_n(k) } \right) - \frac{i}{\sqrt{|\delta|}} e^{-\frac{1}{2}i \arg \delta} \left( \frac{1}{\sqrt{2k_0d}}+ O(\delta)\right) \\
&=:& \gamma_0 - \frac{i}{\sqrt{|\delta|}} e^{-\frac{1}{2}i \arg \delta} \left( \frac{1}{\sqrt{2k_0d}}+ O(\delta)\right).
\end{eqnarray*}
Consequently,  from \eqref{eq-formula_p} it follows that
\begin{eqnarray*}
p(k, \kappa^0, d \varepsilon) &=& \varepsilon +\left[ \dfrac{\cot k }{k} + \dfrac{1}{k\sin k} +  \varepsilon \gamma(k, \kappa^0, d) + \dfrac{1}{\pi} \varepsilon \ln \varepsilon   \right]   \left(\alpha + r(k, \varepsilon) \right), \\
&=& \varepsilon +\left[ \dfrac{\cot k }{k} + \dfrac{1}{k\sin k} +  \varepsilon \gamma_0 - \frac{i \, \varepsilon}{\sqrt{|\delta|}} e^{-\frac{1}{2}i \arg \delta} \left( \frac{1}{\sqrt{2k_0d}}+ O(\delta)\right) + \dfrac{1}{\pi} \varepsilon \ln \varepsilon   \right]   \left(\alpha + r(k, \varepsilon) \right)
\end{eqnarray*}
Similarly, using $\eqref{eq-formula_q}$ for $q(k, \kappa^0, d \varepsilon)$ yields
\begin{equation*}
q(k, \kappa^0, d \varepsilon) = \varepsilon +\left[ \dfrac{\cot k }{k} - \dfrac{1}{k\sin k} +  \varepsilon \gamma_0 - \frac{i \, \varepsilon}{\sqrt{|\delta|}} e^{-\frac{1}{2}i \arg \delta} \left( \frac{1}{\sqrt{2k_0d}}+ O(\delta)\right) + \dfrac{1}{\pi} \varepsilon \ln \varepsilon   \right]   \left(\alpha + s(k, \varepsilon) \right).
\end{equation*}

\begin{lem} \label{lem-pq_Rayleigh_anomaly}
Assume that $\delta=C_0 \cdot \varepsilon^{2\tau}$ wherein $C_0$ is some constant and $0<\tau<1$, and $kd = k^0 d+\delta$.
If $k=m\pi$, then for $\delta>0$,
$$ p(k; \kappa^0, d, \varepsilon)= -\dfrac{i \alpha}{\sqrt{2|C_0|k_0d}} \cdot \varepsilon^{1-\tau}+\min\{O( \varepsilon^{2-2\tau}), O( \varepsilon) \} $$
and
$$ q(k; \kappa^0, d, \varepsilon)= -\dfrac{i \alpha}{\sqrt{2|C_0|k_0d}} \cdot \varepsilon^{1-\tau}+\min\{O( \varepsilon^{2-2\tau}), O( \varepsilon) \} .$$
when $m$ is odd and even respectively. For $\delta<0$,
$$ p(k; \kappa^0, d, \varepsilon)= -\dfrac{\alpha}{\sqrt{2|C_0|k_0d}} \cdot \varepsilon^{1-\tau}+\min\{O( \varepsilon^{2-2\tau}), O( \varepsilon) \} $$
and
$$ q(k; \kappa^0, d, \varepsilon)= -\dfrac{\alpha}{\sqrt{2|C_0|k_0d}} \cdot \varepsilon^{1-\tau}+\min\{O( \varepsilon^{2-2\tau}), O( \varepsilon) \} .$$
when $m$ is odd and even respectively.
\end{lem}
The lemma follows by noting that $\dfrac{\cot k }{k} + \dfrac{1}{k\sin k} =0$ and $\dfrac{\cot k }{k} - \dfrac{1}{k\sin k} =0$ when $k=m\pi$ for odd and even $m$ respectively.

From the discussions in Section 2, the field enhancement will occur for small $p$ or $q$.
If the resonant frequency $m\pi$ is close to a cut-off frequency $k^0$, then from Lemma \ref{lem-pq_Rayleigh_anomaly}
and the same calculation as in Section \ref{sec-field_res}, it can be seen that
the magnitude of the near-field wave is
$O(\varepsilon^{\tau-1})$ and $O(\varepsilon^{\tau}\ln \varepsilon)$ in the silts and on the slit apertures, respectively.
Therefore, in contrast to the case when the resonant frequency is away from the Rayleigh cut-off frequencies,  for which the field is enhanced by an order of $O(\varepsilon^{-1})$, the field enhancement becomes weaker if the resonant frequency is close the Rayleigh cut-off frequencies.
In addition, from the above lemma, it is observed that the wave field at the resonance frequency  has  a phase difference of $\pi$ for $\delta>0$ and $\delta<0$.

\section{A discussion on embedded eigenvalues}\label{sec-emb_eig}
As discussed in Section \ref{sec-eig_res}, for each $\kappa\in(-\pi/d,\pi/d]$, there exist real eigenvalues $k_m$ such that $k_m<|\kappa|$
and is below the light line. The corresponding eigenmodes are surface bound states that are confined near the periodic structure.
In addition, the dispersion curve $k_m(\kappa)$ is continuous (cf. Lemma \ref{lem-Lambda}), and such surface bound states are robust in the sense that
they persist if $\kappa$ is perturbed. For the periodic structure, there may exist real eigenvalues that satisfy
$k>|\kappa|$ and lie above the light line. Such eigenvalues are embedded in the continuous spectrum, and they coincide with the intersection point of 
the complex dispersion relation for the quasi-modes \cite{shipman07, shipman03, shipman10} and the real $(\kappa, k)$-plane. Especially, the corresponding eigenmodes are not robust with respect to perturbation of $\kappa$.
The dissolution of embedded eigenvalues in the continuous spectrum is the mechanism behind transmission anomaly and field enhancement for the periodic slab structure (for instance, Fano resonance) 
when it is illuminated by a plane wave. We refer to \cite{shipman10} and references therein for detailed discussion.

For the periodic structure considered in this paper, unfortunately,  such embedded eigenvalues 
do not exist for the scattering operator $A(\kappa, d, \varepsilon)$ when the domain of the operator is restricted to the space of quasi-periodic functions such that
\begin{equation*}
u_{\varepsilon}(x_1+d,x_2)=e^{i\kappa d} u_{\varepsilon}(x_1,x_2),
\end{equation*}
or more precisely, the function space
$$
H^1_{\kappa,d }(\Omega_{\varepsilon}):= \Big\{u: u \in H^1(\Omega_{\varepsilon}), u(x_1+d, x_2)= e^{i\kappa d} u(x_1, x_2)\Big\}.
$$
Indeed, as discussed in Section \ref{sec-eig_res}, the associated dispersion relation is determined by the roots of the functions $p(k; \kappa, d, \varepsilon)$ or $q(k; \kappa, d, \varepsilon)$, and its asymptotic expansion is given by \eqref{eqn-asym_res_eig}. Moreover, as stated in Section \ref{sec-asy_res_eig},
$\Im\,\gamma(m\pi, \kappa, d)$ holds as long as $m\pi>|\kappa|$. Therefore, $k_m$ will not be a real eigenvalue 
if it lies above the light line or when $k_{m}>|\kappa|$.

However, if we view the periodic structure with a period of $2d$ instead of $d$, and seek for quasi-periodic solutions in $H^1_{\kappa, 2d}(\Omega_{\varepsilon})$ such that
\begin{equation*}
u_{\varepsilon}(x_1+2d,x_2)=e^{i\kappa 2d} u_{\varepsilon}(x_1,x_2),
\end{equation*}
then embedded eigenvalues may exist for the corresponding scattering operator $A(\kappa, 2d, \varepsilon)$.
In what follows, we explore the embedded eigenvalues of the operator $A(\kappa, 2d, \varepsilon)$, which follows the construction in \cite{bonnet_starling94}.

Let us start with an elementary observation. For each fixed $\kappa\in(-\pi/d,\pi/d]$, we define $\hat{\kappa}$ to be one of the two numbers, $\kappa + \pi/d$ and $\kappa - \pi/d$, for which $\hat{\kappa} \in (-\pi/d,\pi/d]$. The Hilbert space $H^1_{\hat{\kappa}, d}(\Omega_{\varepsilon})$ is then well-defined,
and it is clear that $H^1_{\hat{\kappa}, d}(\Omega_{\varepsilon}) \subset H^1_{\kappa, 2d}(\Omega_{\varepsilon})$.

\begin{lem}
The Hilbert space $H^1_{\kappa, 2d}(\Omega_{\varepsilon})$ is the orthogonal sum of the two subspaces $H^1_{\kappa, d}(\Omega_{\varepsilon})$ and $H^1_{\hat{\kappa}, d}(\Omega_{\varepsilon}).$
\end{lem}

\noindent\textbf{Proof}. 
It is clear that $H^1_{\kappa, d}(\Omega_{\varepsilon})$ and $H^1_{\hat{\kappa}, d}(\Omega_{\varepsilon})$ are subspaces of $H^1_{\kappa, 2d},(\Omega_{\varepsilon})$ and they are orthogonal to each other. We only need to show that any functions in $H^1_{\kappa, 2d}(\Omega_{\varepsilon})$ can be written as the sum of two functions in $H^1_{\kappa, d}(\Omega_{\varepsilon})$ and $H^1_{\hat{\kappa}, d}(\Omega_{\varepsilon})$ respectively. Without loss of generality, let us take $\hat{\kappa}=\kappa + \pi/d$. 

Let $f \in H^1_{\kappa, 2d}(\Omega_{\varepsilon})$. Then $f$ adopts the following representation:
\[
f(x_1, x_2) = \sum_{n} f_n(x_2) e^{i\kappa x_1+ i\frac{n\pi}{d}x_1}
\]
when $f_n$'s are expansion coefficients.  By rearranging the terms in the above series, we see that
\[
f 
=\sum_{n} f_{2n}(x_2)e^{i\kappa x+ i\frac{2n\pi}{d}x_1}+\sum_{n} f_{2n+1}(x_2)e^{i\hat{\kappa} x+ i\frac{2n\pi}{d}x_1}.
\]
It is clear that the former belongs to $H^1_{\kappa, d}(\Omega_{\varepsilon})$ and the latter belongs to $H^1_{\hat{\kappa}, d}(\Omega_{\varepsilon})$.
This completes the proof of the lemma. 

As a consequence of the above lemma, we see that the embedded eigenvalues for the operator $A(\kappa, 2d, \varepsilon)$ of the 2d-periodic scattering problem are either associated with the eigenvalues of the operator $A(\kappa, d, \varepsilon)$ or those of the operator $A(\hat\kappa, d, \varepsilon)$.
Equivalently, these correspond to
the roots of the functions $p(k, \kappa, d, \varepsilon)$ and $q(k, \kappa, d, \varepsilon)$, or the functions $p(k, \hat{\kappa}, d, \varepsilon)$ or $q(k, \hat{\kappa}, d, \varepsilon)$. 
We now give a concrete example for such embedded eigenvalues. 
From Remark \ref{rmk-11},  it is known that if
$$ k_1(\pi/d, d, \varepsilon) < \pi/d,$$
then $k_1(\pi/d, d, \varepsilon)$ is an eigenvalue of the operator $A(\kappa, d, \varepsilon)$ for $\kappa=\pi/d$,  and the corresponding surface bound state $\tilde{u}_{\varepsilon}$ satisfies
$$
\tilde{u}_{\varepsilon}(x_1+d,x_2)=e^{i \pi} \tilde{u}_{\varepsilon}(x_1,x_2).
$$
Therefore, it follows
$$
\tilde{u}_{\varepsilon}(x_1+2d,x_2)= \tilde{u}_{\varepsilon}(x_1,x_2).
$$
This shows that $k=k_1(\pi/d, d, \varepsilon)$ is also an eigenvalue of the operator $A(\hat\kappa, 2d, \varepsilon)$ for $\hat\kappa =0$,
and it is an embedded eigenvalue that satisfies $k_1(\pi/d, d, \varepsilon)>|\hat\kappa|=0$.
We may draw similar conclusion for those $(k, \kappa)$ near $(\pi/d, 0)$. 

In general, if we view the periodic structure with a period of $md$, where $m$ is an arbitrary positive integer, and seek for the quasi-periodic solutions in $H^1_{\kappa, md}(\Omega_{\varepsilon})$ such that
\begin{equation*}
u_{\varepsilon}(x_1+md,x_2)=e^{i\kappa md} u_{\varepsilon}(x_1,x_2),
\end{equation*}
a discussion similar to the $2d$-period scattering problem would show the existence of embedded eigenvalues for the corresponding operator
$A(\kappa, md, \varepsilon)$.
However, all these embedded eigenvalues are constructed for the $d$-period scattering problem, and consequently they are robust.
Transmission anomaly and field enhancement would not necessarily occur when the periodic slab structure is illuminated by a plane wave with 
a frequency given by an embedded eigenvalue.
Since, as stated at the beginning of this section, the transmission anomaly such as Fano resonant phenomenon occurs when
the eigenvalue is dissolved into complex-valued resonances, and the corresponding surface bound state is not robust with respect to perturbation of $\kappa$. More precisely,  let $k_0$ be a simple embedded eigenvalue of the operator $A(\kappa_0, md,\varepsilon)$,
then a sufficient condition for Fano resonance to occur is the following (cf. \cite{shipman07, shipman03, shipman10}): 
there exist a small neighbourhood $U\in \mathbf{C}^2$ of $(k_0, \kappa_0)$ such that
$(k, \kappa) = (k_0, \kappa_0)$ is the unique point in $U \cap \mathbf{R}^2$ satisfying 
the dispersion relation
\[
p(k, \hat{\kappa}, d, \varepsilon) =0, \quad \mbox{or}\,\,\,\, 
q(k, \hat{\kappa}, d, \varepsilon) =0.  
\]
That is, $(k, \kappa) = (k_0, \kappa_0)$ is an isolated pair in the real place $\mathbf{R}^2$. However, 
from previous analysis, it is clear that this condition is not satisfied for embedded eigenvalues constructed for the $d$-period scattering problem.

\appendix
\section{Asymptotic expansion of $G_\varepsilon^i(X, Y)$ and $\tilde G_\varepsilon^i(X, Y)$ }
Recall that
\begin{equation}\label{Gi}
G_\varepsilon^i(X, Y) = \dfrac{1 }{\varepsilon} \sum_{m=0}^\infty\left(\sum_{n=0}^\infty c_{mn}\alpha_{mn} \right)  \cos(m\pi X) \cos(m\pi Y).
\end{equation}
Let $\displaystyle{C_m=\sum_{n=0}^\infty c_{mn}\alpha_{mn}}$. Then from the representation of elementary functions by series,
it can be shown that
$$ C_0(k) =  \sum_{n=1}^\infty \dfrac{2}{k^2 - (n\pi)^2} +  \dfrac{1}{k^2} = \dfrac{\cot{k}}{k}, $$
\begin{eqnarray*}
 C_m(k,\varepsilon) &=& \sum_{n=1}^\infty \dfrac{4}{k^2-(m\pi/\varepsilon)^2 -  (n\pi)^2} + \dfrac{2}{k^2- (m\pi/\varepsilon)^2}  \\
                                &=& -\dfrac{2}{\sqrt{(m\pi/\varepsilon)^2-k^2 } } \coth\left( \sqrt{(m\pi/\varepsilon)^2-k^2 } \right) \\
                                &=& - \dfrac{2\varepsilon}{m\pi} - \dfrac{k^2\varepsilon^3}{m^3\pi^3} + O\left(\dfrac{\varepsilon^5}{m^5}\right), \quad m\ge 1. 
 \end{eqnarray*}
Substituting into \eqref{Gi} yields the desired expansion for $G_\varepsilon^i(X, Y)$ given as follows:
\begin{eqnarray*}
G_\varepsilon^i(X, Y) &=&  \dfrac{1 }{\varepsilon} \bigg\{  C_0(k) -  \sum_{m\ge 1} \dfrac{2\varepsilon}{\pi m} \cos(m\pi X) \cos(m\pi Y) 
-  \sum_{m\ge 1} \dfrac{k^2\varepsilon^3}{m^3\pi^3}  \cos(m\pi X) \cos(m\pi Y) \nonumber \\ && +  O\left( \sum_{m\ge 1} \dfrac{\varepsilon^5}{m^5}\right) \bigg\}  \\
&=&  \dfrac{\cot k }{k \varepsilon} + \left(-\dfrac{2}{\pi}\right)\left[ -\ln 2 - \dfrac{1}{2}\ln \left(\abs{\sin \left(\frac{\pi(X+Y)}{2}\right)}\right) - \dfrac{1}{2}\ln \left(\abs{\sin \left(\frac{\pi(X-Y)}{2}\right)}\right) \right]  \nonumber \\
 && + O(k^2\varepsilon^2).
\end{eqnarray*}
On the other hand,
\begin{equation*}
\tilde G_\varepsilon^i(X, Y) = \dfrac{1 }{\varepsilon} \sum_{m=0}^\infty\left(\sum_{n=0}^\infty (-1)^n c_{mn}\alpha_{mn} \right)  \cos(m\pi X) \cos(m\pi Y).
\end{equation*}
Let $\displaystyle{\tilde C_m=\sum_{n=0}^\infty (-1)^n c_{mn}\alpha_{mn}}$. Again from the following representation of elementary functions as series
$$ \tilde C_0(k) =  \sum_{n=1}^\infty \dfrac{2(-1)^n}{k^2 - (n\pi)^2} +  \dfrac{1}{k^2} = \dfrac{1}{k\sin{k}}, $$
\begin{eqnarray*}
\tilde C_m(k,\varepsilon) &=& \sum_{n=1}^\infty \dfrac{(-1)^n \cdot 4}{k^2-(m\pi/\varepsilon)^2 -  (n\pi)^2} + \dfrac{2}{(m\pi/\varepsilon)^2 }  \\
                                &=& -\dfrac{2}{\sqrt{(m\pi/\varepsilon)^2-k^2 } \sinh\left( \sqrt{(m\pi/\varepsilon)^2-k^2 } \right)} \\
                                &=&  O\left(\dfrac{\varepsilon}{m\pi}e^{-m\pi/\varepsilon}\right),  \quad \quad m\ge1,
\end{eqnarray*}
we obtain
\begin{equation*}
\tilde G_\varepsilon^i(X, Y) =  \dfrac{1}{(k\sin k) \varepsilon} + O\left(e^{-1/\varepsilon}\right).
\end{equation*}

\bibliography{references}

\end{document}